\documentclass[11pt,a4paper]{article}

\usepackage[T1]{fontenc}
\usepackage[utf8]{inputenc}
\usepackage{lmodern}
\usepackage{amsmath,amssymb,amsthm,mathtools}
\usepackage{booktabs}
\usepackage{enumitem}
\usepackage{geometry}
\usepackage{xcolor}
\usepackage{listings}
\usepackage[hidelinks]{hyperref}
\usepackage[nameinlink,capitalise,noabbrev]{cleveref}
\usepackage{algorithm}
\usepackage{algpseudocode}
\usepackage{float}
\usepackage{tikz}
\usepackage{soul}
\usepackage[most]{tcolorbox}
\usepackage{subcaption}

\newtheorem{theorem}{Theorem}[section]
\newtheorem{lemma}[theorem]{Lemma}
\newtheorem{proposition}[theorem]{Proposition}
\newtheorem{corollary}[theorem]{Corollary}
\theoremstyle{definition}

\newtheorem{example}[theorem]{Example}
\theoremstyle{remark}

\newtheorem{conjecture}[theorem]{Conjecture}

\DeclareMathOperator{\diam}{diam}
\newcommand{\F}{F^s_{k\times 1}}

\title{On generalized $k$-token graphs with distinguishable tokens and a maximum capacity of $s$ tokens per vertex}
\author{Yiting Cai\\
\small{School of Mathematical Sciences, South China Normal University, Guangzhou, China}\\
\small{Departament de Matem\`atica, Universitat de Lleida, Igualada (Barcelona), Catalonia}\\
\vspace{.2cm}
\small{\texttt{yitingcai@m.scnu.edu.cn}}\\
Cristina Dalf\'o\\
\small{Departament de Matem\`atica, Universitat de Lleida, Igualada (Barcelona), Catalonia}\\
\vspace{.2cm}
\small{\texttt{cristina.dalfo@udl.cat}}\\
Blas Fern\'andez\\
\small{University of Primorska, Koper, Slovenia}\\
\vspace{.2cm}
\small{\url{blas.fernandez@famnit.upr.si}}\\
Miquel \`Angel Fiol\\
\small{Departament de Matem\`atiques, Universitat Polit\`ecnica de Catalunya, Barcelona, Catalonia;}\\
\small{Barcelona Graduate School of Mathematics, Barcelona, Catalonia;}\\
\small{Institut de Matem\`atiques de la UPC-BarcelonaTech (IMTech), Barcelona, Catalonia}\\
\small{\texttt{miguel.angel.fiol@upc.edu}}
}
\date{}

\begin{document}

\maketitle

\begin{abstract}
Supertoken graphs \(F^s_{k\times1}(G)\), also known as generalized token graphs, model configurations of \(k\) distinguishable tokens on a graph \(G\), with at most \(s\) tokens at each vertex. We study structural and metric properties of these graphs in terms of the base graph. For regular \(G\), we characterize when \(F^s_{k\times1}(G)\) is regular. We show that bipartiteness is inherited from \(G\) and characterize it for complete base graphs. For \(s\ge2\), we characterize connectedness and show that \(t\)-vertex-connectivity of \(G\) implies \(t\)-vertex-connectivity of the supertoken graph. We obtain exact diameter formulas in several regimes, including complete graphs and stars, and derive lower bounds and a conjectured formula for odd cycles. We also determine the clique number, give bounds and equality conditions for the chromatic number, and characterize planarity for complete base graphs.

\end{abstract}




\section{Introduction and definitions}
\label{sec:intro}

Token graphs provide a natural configuration model for the motion of
several tokens on a graph. Given a base graph $G$, each configuration
records the positions of the tokens, and two configurations are adjacent
when one token can be moved along an edge of $G$. Thus, the resulting
graph reflects both the structure of $G$ and the interaction among the
tokens.

Classical token graphs allow at most one token at each vertex.
Supertoken graphs replace this exclusion rule by a capacity constraint,
allowing up to $s$ tokens per vertex. When the tokens are distinguishable,
a configuration also records which tokens occupy each vertex. 

Capacity constraints introduce new phenomena. A token may be unable to
follow a shortest path because its next vertex is full, forcing other
tokens to move temporarily to create a vacant place. Consequently, even
simple base graphs can produce supertoken graphs with nontrivial metric
and structural behavior.

We aim to study how the base graph and the parameters $k$ and $s$
determine fundamental properties of the resulting supertoken graph. We
focus on bipartiteness, connectivity, distances and diameter, clique and chromatic number, and planarity, with
particular attention to complete graphs, stars, and cycles.

Supertoken graphs (also called generalized token graphs) generalize token graphs, which were previously known as symmetric powers (see Audenaert, Godsil, Royle, and Rudolph \cite{agrr07}). Later, Fabila-Monroy, Flores-Pe\~{n}aloza, Huemer, Hurtado, Urrutia, and Wood~\cite{ffhhuw12}.
 referred to these graphs as \(k\)-token graphs of \(G\).

Here, we consider distinguishable tokens and a common capacity $s$ at
each vertex.


Let \(G=(V,E)\) be a graph of order \(n\), and let \(k,s\ge1\) with \(k\le ns\). A configuration (or vertex of the supertoken graph) is an ordered \(k\)-tuple
\[
\boldsymbol{x}=(x_1,\ldots,x_k)\in V^k
\]
such that no vertex of \(G\) occurs more than \(s\) times. The \(i\)-th coordinate records the position of token \(T_i\). Two configurations \(\boldsymbol{x}=(x_1,\ldots,x_k)\) and \(\boldsymbol{y}=(y_1,\ldots,y_k)\) are adjacent if there is a unique \(i\) such that \(x_i\ne y_i\), \(x_i y_i\in E(G)\), and \(x_j=y_j\) for every \(j\ne i\). 

For a configuration $\boldsymbol{x}=(x_1,\ldots,x_k)$, let $a_W(v)$ denote the number of tokens at $v\in V(G)$, that is 
\[
a_{\boldsymbol{x}}(v):=|\{i: x_i=v\}|. 
\]
Then $0\le a_{\boldsymbol{x}}(v)\le s$ and 
\[
\sum_{v\in V(G)}a_{\boldsymbol{x}}(v)=k. 
\]
We say that vertex $v$ is full in ${\boldsymbol{x}}$ if $a_{\boldsymbol{x}}(v)=s$. 

The order and size of these supertoken graphs were determined by 
Song, Dalf\'o, Fiol, Mora, and Zhang \cite{sdfmz26}.


An example of $F_{3\times1}^2(P_3)$ is shown in Figure 
\ref{fig:F_{3x1}^2(P_3)}.

\begin{figure}[t]
\begin{center}
\includegraphics*[width=10cm]{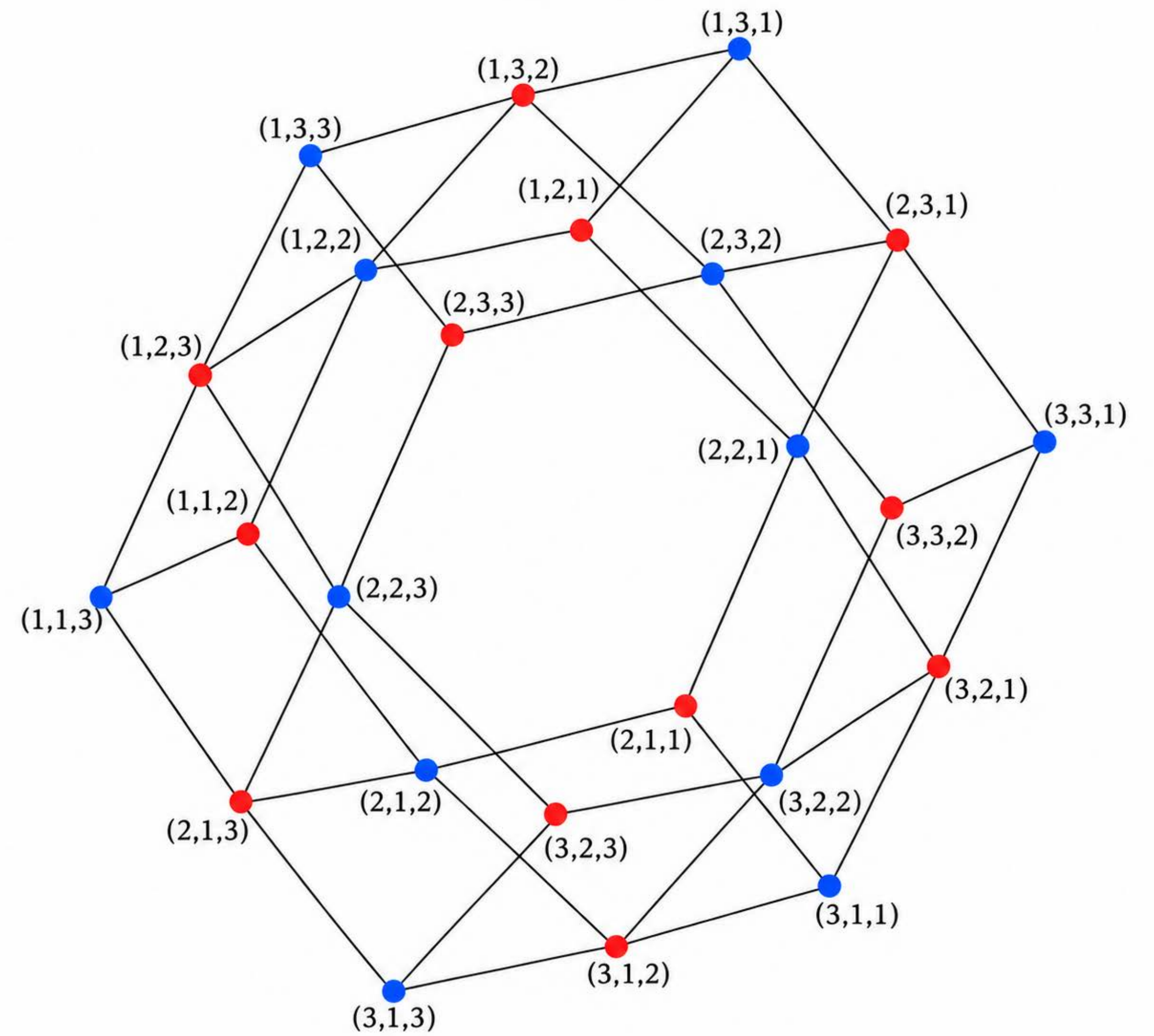}  
\caption{The supertoken graph $F_{3\times1}^2(P_3)$ as a bipartite graph.} 
\label{fig:F_{3x1}^2(P_3)}
\end{center}
\end{figure}

This paper is organized as follows. We investigate how structural properties of \(F^s_{k\times 1}(G)\) are related to those of the underlying base graph \(G\). In Section \ref{sec:basic-prop}, for regular base graphs, we characterize precisely when the corresponding supertoken graph is regular. Besides, we show that bipartiteness is inherited from the base graph. In Section \ref{sec:distances}, we obtain exact diameter formulas in several regimes, including complete graphs and stars, and derive lower bounds and a conjectured exact formula for odd cycles. Sections \ref{sec:connect}, \ref{sec:clique}, and \ref{sec:chromatic} are devoted, respectively, to vertex-connectivity, clique number, and chromatic number. In particular, we 
prove that vertex-connectivity is inherited from the base graph (in the sense that \(t\)-connectivity of \(G\) implies \(t\)-connectivity of the supertoken graph), and establish bounds and characterizations for the clique and chromatic numbers. Finally, in Section \ref{sec:planar}, we characterize planarity when the base graph is complete. The Appendix contains the pseudo-code for computing the diameter of \(F^s_{k\times 1}(G)\).



\section{Basic structural properties}
\label{sec:basic-prop}

\subsection{Degree and regularity}


Regularity is a natural structural property to investigate for
supertoken graphs. Even when the base graph is regular, however, the
capacity constraint may create configurations with different numbers
of available token moves. The purpose of this section is to isolate
precisely when this phenomenon occurs.

We begin with a degree formula valid for an arbitrary base graph.
It separates the contribution coming from the degrees in the base
graph from the contribution caused by full vertices, which block
otherwise possible moves. We then specialize to regular base graphs.
In that setting, the unrestricted contribution is constant, and all
variation in the degree is controlled by the capacity constraint.

Let $G$ be a simple, undirected graph on $n$ vertices. Let $s\geq2$, and 
$1\leq k\leq ns$. Label the $k$ distinguishable tokens by
$T_1,\ldots,T_k$. A vertex of $F^s_{k\times1}(G)$ is an admissible
configuration $\boldsymbol{x}=(x_1,\ldots,x_k)$,
where $x_i\in V(G)$ is the vertex occupied by token $T_i$.
Several tokens may occupy the same vertex of $G$, but no vertex
may contain more than $s$ tokens.

For a configuration $\boldsymbol{x}=(x_1,\ldots,x_k)$ and $v\in V(G)$, define
\[
S(\boldsymbol{x}):=\{v\in V(G):a_{\boldsymbol{x}}(v)=s\}.
\]


Accordingly, a token at $u$ may move to a neighbor $v$ of $u$
precisely when $v\notin S(\boldsymbol{x})$, since its destination
must have room for one additional token.

\subsubsection{A degree formula}

The following lemma provides the basic counting tool for this section.
The first expression counts, for each occupied base vertex, the
destinations that are still available. The second expression presents
the same count from a complementary viewpoint: start with all possible
moves and subtract those whose destination is full.

\begin{lemma}
\label{lem:supertoken-degree}
Let $G$ be a finite simple graph, let $s\geq1$, and let
$1\leq k\leq s|V(G)|$. For every vertex $\boldsymbol{x}$ of
$F^s_{k\times1}(G)$,
\[
\deg_{F^s_{k\times1}(G)}(\boldsymbol{x})
=
\sum_{u\in V(G)}
a_{\boldsymbol{x}}(u)
\bigl(\deg_G(u)-|N_G(u)\cap S(\boldsymbol{x})|\bigr).
\]
Equivalently,
\[
\deg_{F^s_{k\times1}(G)}(\boldsymbol{x})
=
\sum_{u\in V(G)}
a_{\boldsymbol{x}}(u)\deg_G(u)
-
\sum_{v\in S(\boldsymbol{x})}
\sum_{u\in N_G(v)}a_{\boldsymbol{x}}(u).
\]
\end{lemma}

\begin{proof}
Each of the \(a_{\boldsymbol{x}}(u)\) tokens at \(u\) can be moved to precisely
$\deg_G(u)-|N_G(u)\cap S(\boldsymbol{x})|$ neighbors of \(u\). Hence, 
\[
\deg_{F^s_{k\times1}(G)}(\boldsymbol{x})=
\sum_{u\in V(G)}
a_{\boldsymbol{x}}(u)\bigl(\deg_G(u)-|N_G(u)\cap S(\boldsymbol{x})|\bigr).
\]
Finally,
\[
\sum_u a_{\boldsymbol{x}}(u)|N_G(u)\cap S(\boldsymbol{x})|
=
\sum_{v\in S(\boldsymbol{x})}\sum_{u\in N_G(v)}a_{\boldsymbol{x}}(u),
\]
by reversing the order of summation.
\end{proof}

The second expression becomes especially transparent when the base
graph is regular: the total number of potential moves, before imposing
the capacity restriction, is $kd$. Thus, the only source of
degree variation is the number of moves blocked by full vertices.




\begin{lemma}\label{cor:blocking}
Let $G$ be a finite simple $d$-regular graph, and let
$H=F^s_{k\times1}(G)$. For every $\boldsymbol{x}\in V(H)$, define
\[
B(\boldsymbol{x}):=\sum_{v\in S(\boldsymbol{x})}\sum_{u\in N_G(v)}a_{\boldsymbol{x}}(u).
\]
Then
\[
\deg_H(\boldsymbol{x})=kd-B(\boldsymbol{x}).
\]
Consequently, $H$ is regular if and only if $B(\boldsymbol{x})$ is independent of
$\boldsymbol{x}\in V(H)$. In particular, if $B(\boldsymbol{x})=b$ for every $\boldsymbol{x}\in V(H)$, then
$H$ is $(kd-b)$-regular.
\end{lemma}

\begin{proof}
The degree formula follows immediately from Lemma~\ref{lem:supertoken-degree}, since
$G$ is $d$-regular and
\[
\sum_{u\in V(G)}a_{\boldsymbol{x}}(u)=k.
\]
The remaining statements are immediate.
\end{proof}


We call $B(x)$ the \emph{blocking number} of $x$; it counts the moves
prohibited because their destination is full. Thus, for a regular base
graph, regularity of the supertoken graph is equivalent to the blocking
number being independent of the configuration.

Therefore, it remains to show that,
except in a few
extremal regimes, one can construct configurations with different
blocking numbers. 
There are two qualitatively distinct situations.
If sufficiently many positions are available, one may compare a
configuration with no full vertex to one containing a full vertex.
Near total capacity, every configuration necessarily has some full
vertices; in that case, the distribution of the vacancies becomes the
relevant parameter.

The next lemma handles this latter regime. 
 Its role is worth noting:
it does not merely show that the set of full vertices can vary, but
that this variation changes the actual number of blocked moves. This
is precisely what is needed for the degree argument.

\begin{lemma}
\label{lem:vacancy-variation}
Let \(G\) be a finite simple \(d\)-regular graph on \(n\geq2\)
vertices, where \(d\geq1\), and let \(s\geq2\). For every
\(h\in\{2,\ldots,n-1\}\), there exist two admissible configurations
with \(ns-h\) tokens having distinct blocking numbers.
\end{lemma}
\begin{proof}
Put $
m:=h-2.
$
We first consider the case in which either $s\geq3$, or $s=2$ and
$d\geq2$. Choose an edge $uv$ of $G$. If $s=2$ and $m\geq d-1$, choose a neighbor
$
u_0\in N_G(u)\setminus\{v\},
$
which exists because $d\geq2$, and choose
$
W\subseteq V(G)\setminus\{u,v\}$ with $|W|=m$
such that $u_0\in W$. This is possible because $m\leq n-3$. In all
other cases, choose any set
$
W\subseteq V(G)\setminus\{u,v\}$ with $|W|=m$. Define configurations $\boldsymbol{x}$ and $\boldsymbol{y}$ by
$$
a_{\boldsymbol{x}}(z)
=
\begin{cases}
s-1, & z\in W\cup\{u,v\},\\
s, & z\notin W\cup\{u,v\},
\end{cases}
\qquad
a_{\boldsymbol{y}}(z)
=
\begin{cases}
s-2, & z=u,\\
s-1, & z\in W,\\
s, & z\notin W\cup\{u\}.
\end{cases}
$$
The configuration $\boldsymbol{x}$ has one vacant position at every
vertex of $W\cup\{u,v\}$, whereas $\boldsymbol{y}$ has two vacant
positions at $u$ and one vacant position at every vertex of $W$.
Thus, both configurations have 
$
|W|+2=m+2=h
$
vacant positions, and hence both have $ns-h$ tokens.

Put
$
T:=V(G)\setminus\bigl(W\cup\{u,v\}\bigr)$, $A:=|N_G(u)\cap W|$, and $C:=|N_G(v)\cap W|.
$ Then $
S(\boldsymbol{x})=T$ and $S(\boldsymbol{y})=T\cup\{v\}.
$
Since the occupancies of $\boldsymbol{x}$ and $\boldsymbol{y}$ differ
only at $u$ and $v$, we have
$
a_{\boldsymbol{y}}(u)-a_{\boldsymbol{x}}(u)=-1$ and
$
a_{\boldsymbol{y}}(v)-a_{\boldsymbol{x}}(v)=1.
$
Moreover,
$
|N_G(u)\cap T|=d-1-A$ and
$
|N_G(v)\cap T|=d-1-C.
$
It follows that
$$
\begin{aligned}
\sum_{z\in T}\sum_{w\in N_G(z)}
\bigl(a_{\boldsymbol{y}}(w)-a_{\boldsymbol{x}}(w)\bigr)
&=-|N_G(u)\cap T|+|N_G(v)\cap T|\\
&=
-(d-1-A)+(d-1-C)=A-C.
\end{aligned}
$$
Furthermore, $v$ is full in $\boldsymbol{y}$ but not in
$\boldsymbol{x}$. Its neighbors consist of $u$, the $C$ vertices in
$N_G(v)\cap W$, and $d-1-C$ further full vertices. Therefore,
$$
\begin{aligned}
\sum_{w\in N_G(v)}a_{\boldsymbol{y}}(w)
&=
(s-2)+C(s-1)+(d-1-C)s=
sd-2-C.
\end{aligned}
$$
Consequently, $B(\boldsymbol{y})-B(\boldsymbol{x})
=
(A-C)+(sd-2-C)
=
sd-2+A-2C.$ Since $u\notin W$, we have $C\leq d-1$. Hence
$$
B(\boldsymbol{y})-B(\boldsymbol{x})
\geq
sd-2-2(d-1)
=
d(s-2).
$$
Thus, if $s\geq3$, then $
B(\boldsymbol{y})-B(\boldsymbol{x})>0.
$ It remains to consider the case $s=2$ and $d\geq2$. If $m<d-1$, then
$C\leq |W|=m$, and hence
$
B(\boldsymbol{y})-B(\boldsymbol{x})
\geq
2d-2-2m
>0.
$
If $m\geq d-1$, then our choice of $W$ gives $A\geq1$, while
$C\leq d-1$. Therefore,
$
B(\boldsymbol{y})-B(\boldsymbol{x})
\geq
2d-2+1-2(d-1)
=
1.
$
Thus, in all cases with $d\geq2$, the configurations
$\boldsymbol{x}$ and $\boldsymbol{y}$ have distinct blocking numbers.

Finally, suppose that $s=2$ and $d=1$. Then $G$ is a disjoint union of
copies of $K_2$, and hence $n$ is even. Since $h\geq2$ and
$h\leq n-1$, we have $n\geq4$. Choose non-adjacent vertices $u$ and
$v$, and let $u'$ and $v'$ be their respective unique neighbors.
Choose $
W\subseteq V(G)\setminus\{u,v\}$ with $|W|=h-2$
such that $v'\notin W$. This is possible because $h-2\leq n-3$.
Define $\boldsymbol{x}$ and $\boldsymbol{y}$ as above. The vertex $v$ is full in $\boldsymbol{y}$ and contributes $2$ to
$B(\boldsymbol{y})$, while the contribution of $v'$ increases by $1$.
If $u'\in W$, then $u'$ is non-full in both configurations, and
therefore
$
B(\boldsymbol{y})-B(\boldsymbol{x})=3.
$
If $u'\notin W$, then $u'$ is full in both configurations and its
contribution decreases by $1$, so
$
B(\boldsymbol{y})-B(\boldsymbol{x})=2.
$
Thus, in either case,
$
B(\boldsymbol{y})-B(\boldsymbol{x})\geq2$ and so, the blocking numbers are distinct.

Since the tokens are distinguishable, each occupancy vector above can be
realized by distributing the $ns-h$ tokens among the vertices according
to the prescribed occupancies. This completes the proof.
\end{proof}

\subsubsection{Classification for regular base graphs}

We can now give the complete classification. The three regular cases
have a simple interpretation. 
When $k\le s$, no full vertex blocks a move, so every configuration
has degree $kd$.
When
$k=ns-1$, there is a unique vacant position, and every permitted move
is directed into that position. Finally, when $k=ns$, every vertex is
full and no move is possible.

All remaining values of $k$ allow enough freedom to produce
configurations with different blocking numbers, and hence different
degrees.

\begin{theorem}
\label{thm:regularity-regular-base}
Let $G$ be a finite simple $d$-regular graph on $n\geq2$ vertices,
where $d\geq1$. Let $s\geq2$ and $1\leq k\leq ns$. 
Then
$F^s_{k\times1}(G)$ is regular if and only if
$
k\leq s$ or 
$
k\in\{ns-1,ns\}.
$
In these cases,  its valency is respectively $kd$, $sd$, and $0$.
\end{theorem}

\begin{proof}
For \(\boldsymbol{x}\in V(F^s_{k\times1}(G))\), let $B(\boldsymbol{x})$ be the blocking number of the
configuration $\boldsymbol{x}$. Suppose first that \(1\leq k\leq s\). 
If \(k<s\), then \(S(x)=\varnothing\). If \(k=s\) and \(v\in S(x)\), then all \(k=s\) tokens are located at \(v\); hence every neighbor of \(v\) is empty and
\[
\sum_{u\in N_G(v)}a_x(u)=0.
\]
Thus \(B(x)=0\) in either case, and consequently every configuration has degree \(kd\).
Now suppose that \(k=ns-1\). Since the total capacity of the \(n\)
vertices of \(G\) is \(ns\), every configuration has exactly one
vacant position. Hence there is a unique vertex \(v\) with
$
a_{\boldsymbol{x}}(v)=s-1,
$
whereas every other vertex is full. Therefore
$
S(\boldsymbol{x})=V(G)\setminus\{v\}.
$
For every full vertex $w$, all its neighbors are full except
possibly $v$. Equivalently, the only possible token moves are moves
from a neighbor of $v$ into $v$. Since $v$ has $d$ neighbors
and each of them contains $s$ tokens, there are exactly $sd$ such
moves. Hence $\deg_{F^s_{k\times1}(G)}(\boldsymbol{x})=sd$
for every $\boldsymbol{x}$. If $k=ns$, every vertex is full, so no token can be moved. Thus $F^s_{k\times1}(G)$ is $0$-regular.

It remains to prove that no other value of \(k\) gives a regular
supertoken graph. Assume that $s<k\leq ns-2$. First suppose that $k\leq n(s-1)$. Choose a configuration
$\boldsymbol{x}$ in which no vertex is full; this is possible
because the $k$ tokens can be distributed with at most $s-1$
tokens at each vertex. Then $S(\boldsymbol{x})=\varnothing$ and hence $B(\boldsymbol{x})=0$. Since $d\geq1$, choose an edge $uv$ of $G$. There is also a configuration $\boldsymbol{y}$ with
$
a_{\boldsymbol{y}}(u)=s$ and $a_{\boldsymbol{y}}(v)\geq1$.
Indeed, place $s$ tokens at $u$, one token at $v$, and distribute
the remaining $k-s-1$ tokens arbitrarily among the remaining
available positions. Since $u\in S(\boldsymbol{y})$ and $v\in N_G(u)$,
we obtain
$
B(\boldsymbol{y})
\geq a_{\boldsymbol{y}}(v)
\geq1.
$
Thus, \(B(\boldsymbol{x})\ne B(\boldsymbol{y})\), and Lemma~\ref{cor:blocking} shows that the supertoken graph is
not regular. Finally, suppose that $n(s-1)<k\leq ns-2$, and put $h:=ns-k$. The inequality $k>n(s-1)$ gives $h<n$, whereas
$k\leq ns-2$ gives $h\geq2$. Hence, $2\leq h\leq n-1$. In other words, every configuration has $h$ vacant positions, where
there are at least two vacancies but not enough to leave every base vertex non-full. By Lemma~\ref{lem:vacancy-variation}, there are configurations
$\boldsymbol{x}$ and $\boldsymbol{y}$ with $ns-h=k$ tokens such
that $B(\boldsymbol{x})\ne B(\boldsymbol{y})$. Thus, by Lemma~\ref{cor:blocking}, the supertoken graph \(F^s_{k\times1}(G)\) is not regular. 

Therefore, \(F^s_{k\times1}(G)\) is regular precisely when
$k\leq s$ or $k\in\{ns-1,ns\}$, with valencies $kd$, $sd$, and $0$, respectively.
\end{proof}

The theorem applies immediately to the standard regular families
below. 

In particular, Theorem \ref{thm:regularity-regular-base} applies to \(K_n\), \(K_{n,n}\), and \(C_n\), giving valencies $k(n-1)$, $kn$, $2k$, respectively, in the range \(k\le s\), with the corresponding values $s(n-1)$, $sn$, $2s$ 
when there is one vacancy.








The assumption \(s\ge2\) is essential: for instance,
\(F^1_{2\times1}(K_4)\) is \(4\)-regular, although \(k=2\) satisfies none of the conditions in Theorem~\ref{thm:regularity-regular-base}.

\subsection{Bipartiteness}
\label{sec:bipartite}

We first show that the bipartiteness of the base graph $G$ is inherited
by $F^s_{k\times1}(G)$ for every admissible choice of $k$ and $s$.
The key observation is that each legal move changes the parity
of the number of tokens occupying one fixed class of a bipartition
of $G$. This yields a bipartition of the graph,
without requiring $G$ or $F^s_{k\times1}(G)$ to be connected.
\begin{theorem}\label{thm:bip}
If $G$ is bipartite, then $\F(G)$ is bipartite for every
admissible $k,s$. This includes disconnected base graphs
and isolated vertices.
\end{theorem}
\begin{proof}
Let $A$ and $B$ be the two classes of a bipartition of $G$,
so that $V(G)=A\cup B$ and $A\cap B=\varnothing$.
Define the parity function
$$
p_A:V(F^s_{k\times1}(G))\longrightarrow\{0,1\}
$$
by
$$
p_A({\boldsymbol{x}}):=
|\{i\in\{1,\ldots,k\}:x_i\in A\}|\pmod 2,
\qquad {\boldsymbol{x}}=(x_1,\ldots,x_k).
$$
Thus, $p_A({\boldsymbol{x}})$ records whether the number of tokens at
vertices of $A$ is even or odd.\\
Every legal move transfers exactly one token from $A$ to $B$
or from $B$ to $A$. Hence, whenever ${\boldsymbol{x}}$ and ${\boldsymbol{y}}$ are adjacent,
$$
p_A({\boldsymbol{y}})=1-p_A({\boldsymbol{x}}).
$$
Therefore, the sets
$$
\mathcal V_0=\{{\boldsymbol{x}}:p_A({\boldsymbol{x}})=0\},
\qquad
\mathcal V_1=\{{\boldsymbol{x}}:p_A({\boldsymbol{x}})=1\}
$$
form a bipartition of $F^s_{k\times1}(G)$.
The argument also applies when $G$ is disconnected or
the graph contains isolated vertices.
\end{proof}
(The converse fails. Although $K_3$ is not bipartite,
$F^1_{2\times1}(K_3)$ is the cycle
$(1,2)$, $(3,2)$, $(3,1)$, $(2,1)$, $(2,3)$, $(1,3)$, $(1,2)$
and, hence, is bipartite. 
Thus, the forward implication holds for every bipartite base graph.)
    
We next characterize when the supertoken graph  $F^s_{k\times1}(K_n)$ is bipartite.
Recall that $K_n$ has $n$ vertices, each
vertex can contain at most $s$ tokens, and the total number of labeled
tokens is $k$. Consequently, the natural parameter range is
$1\leq k\leq ns$. The upper bound $k\leq ns$ follows from the fact that the total
capacity of the $n$ vertices is $ns$.

{When $s=1$, Dalf\'o, Fern\'andez, Fiol, and
Reyes~\cite{dffr26} proved that $F^1_{k\times1}(K_n)$ is bipartite if
and only if $k\in\{n-1,n\}$. 
The situation is different when 
$s\geq2$. As we show below, if $n\geq3$, then every non-full
supertoken graph $F^s_{k\times1}(K_n)$ contains an odd cycle and is
therefore nonbipartite.}

We first consider the exceptional base graph $K_2$. Since $K_2$ is
bipartite, Theorem~\ref{thm:bip} implies that
$F^s_{k\times1}(K_2)$ is bipartite for every admissible choice of
$s$ and $k$.

 We now give the characterization for complete graphs of order at
least three.

\begin{theorem}
	\label{thm:bipartite-complete-supertoken}
	Let $n\geq3$, let $s\geq2$, and let $1\leq k\leq ns$.
	Then $F^s_{k\times1}(K_n)$ is bipartite if and only if $k=ns$.
\end{theorem}

\begin{proof}
	Suppose first that $k=ns$. In every configuration, each of the $n$ vertices of $K_n$ contains
	exactly $s$ tokens. Hence, every vertex of the base graph is full, and
	no token can be moved. Therefore, $F^s_{ns\times1}(K_n)$
	is edgeless and is, consequently, bipartite.

	For the converse, suppose that $1\leq k<ns$.
	We prove that $F^s_{k\times1}(K_n)$ contains an odd cycle. We
	distinguish two cases according to the number of vacant places.
	
	\medskip
	
	\noindent
	\textbf{Case 1:} $k\leq ns-2$.
	
	Choose three distinct vertices $u,v,w\in V(K_n)$, and distinguish one
	of the tokens, say $T_1$. We first place the remaining
	$k-1$ tokens in such a way that each of $u,v,w$ contains at most
	$s-1$ tokens.
		
	To see that such a placement is possible, reserve one place at each
	of the three vertices $u,v,w$. The number of places that remain
	available is
	$
	3(s-1)+(n-3)s=ns-3.
	$
	On the other hand, the assumption $k\leq ns-2$ gives
	$
	k-1\leq ns-3.
	$
	Thus, the tokens $T_2,\ldots,T_k$ can be distributed among the
	vertices of $K_n$ so that no vertex contains more than $s$ tokens
	and, in addition, each of $u,v,w$ contains at most $s-1$ tokens.
	
	Keep the positions of $T_2,\ldots,T_k$ fixed. Let ${\boldsymbol{x}}_u$, ${\boldsymbol{x}}_v$, and
	${\boldsymbol{x}}_w$ be the configurations obtained by placing $T_1$ at $u$, $v$,
	and $w$, respectively.
	
	Since $K_n$ is complete, the vertices $u,v,w$ are pairwise adjacent.
	Moreover, each of them has a place available for $T_1$.
		Furthermore, the three configurations are pairwise adjacent. To pass
	from ${\boldsymbol{x}}_u$ to ${\boldsymbol{x}}_v$, we move only $T_1$ from $u$ to $v$. This move
	is legal because $u$ and $v$ are adjacent in the complete graph
	$K_n$ and, in ${\boldsymbol{x}}_u$, the vertex $v$ contains at most $s-1$ tokens.
	Similarly, moving $T_1$ from $v$ to $w$ gives the edge
	$
	{\boldsymbol{x}}_v{\boldsymbol{x}}_w,
	$
	and moving $T_1$ from $w$ to $u$ gives the edge
	$
	{\boldsymbol{x}}_w{\boldsymbol{x}}_u.
	$
	
	 Consequently,
	the moves
	\[
	{\boldsymbol{x}}_u\longrightarrow {\boldsymbol{x}}_v,\qquad
	{\boldsymbol{x}}_v\longrightarrow {\boldsymbol{x}}_w,\qquad
	{\boldsymbol{x}}_w\longrightarrow {\boldsymbol{x}}_u
	\]
	are all legal. Hence ${\boldsymbol{x}}_u,{\boldsymbol{x}}_v,{\boldsymbol{x}}_w,{\boldsymbol{x}}_u$ forms a $3$-cycle in $F^s_{k\times1}(K_n)$. It follows that
	$F^s_{k\times1}(K_n)$ is not bipartite.
	
	\medskip
	
	\noindent
	\textbf{Case 2:} $k=ns-1$.
	
	Choose three distinct vertices $u,v,w\in V(K_n)$. 	At each of these vertices, place $s-2$ tokens that will remain fixed throughout the construction. At each of the remaining $n-3$ vertices, place $s$ fixed tokens. The number of fixed tokens used in this way is
	$
	3(s-2)+(n-3)s=ns-6.
	$
	Since $k=ns-1$, there remain exactly
	$
	(ns-1)-(ns-6)=5
	$
	tokens. Denote them by
	$
	T_1,T_2,T_3,T_4,T_5.
	$
	The fixed background tokens leave two available places at each of
	$u,v,w$. We place $T_4$ at $v$ and $T_5$ at $w$ and keep them fixed.
	Only the positions of $T_1,T_2,T_3$ will change.
	
	For $a,b,c\in\{u,v,w\}$, write 	$[a,b,c]$
	for the configuration in which $T_1,T_2,T_3$ occupy $a,b,c$,
	respectively, while $T_4$ occupies $v$, $T_5$ occupies $w$, and all
	background tokens remain at their fixed vertices. 	
After fixing \(s-2\) tokens at each of \(u,v,w\), \(s\) tokens elsewhere, and fixing \(T_4\) at \(v\) and \(T_5\) at \(w\), the configurations
\[
[u,u,v], [u,w,v], [u,w,u], [v,w,u],
[v,u,u], [w,u,u], [w,u,v]
\]form a \(7\)-cycle.



Hence,
$F^s_{k\times1}(K_n)$ contains an odd cycle and is therefore not
bipartite.
	
	We have proved that $F^s_{k\times1}(K_n)$ contains an odd cycle for
	every $k<ns$, whereas it is edgeless when $k=ns$. This proves the
	result.
\end{proof}

Combining the preceding theorem with the known result for $s=1$ gives the following complete characterization for $n\geq3$.

\begin{corollary}
	\label{cor:bipartite-complete-all-capacities}
	Let $n\geq3$, let $s\geq1$, and let $1\leq k\leq ns$. Then
	$F^s_{k\times1}(K_n)$ is bipartite if and only if one of the
	following conditions holds:
	\begin{enumerate}[label=\textnormal{\roman*)}]
		\item $s=1$ and $k\in\{n-1,n\}$;
		\item $s\geq2$ and $k=ns$.
	\end{enumerate}
\end{corollary}


\section{Connectivity}
\label{sec:connect}

We next study how connectivity properties of the base graph are reflected in the corresponding supertoken graph. We first characterize when $F^s_{k\times1}(G)$ is connected, and then show that, for $s\ge2$, vertex-connectivity is inherited from the base graph: if $G$ is $t$-vertex-connected, then so is $F^s_{k\times1}(G)$.

\begin{theorem}\label{th:connected}
Let $G$ be a graph of order $n\ge2$, and let $s,k$ be integers with
$2\le s<k\le ns$. Then
$F^s_{k\times1}(G)$
is connected if and only if $G$ is connected and $k<ns$.
\end{theorem}

\begin{proof}
Write $H=F^s_{k\times1}(G)$.

Suppose first that $H$ is connected. If $G$ is disconnected, choose vertices
$u$ and $v$ in distinct components of $G$, and configurations ${\boldsymbol{x}},{\boldsymbol{y}}\in V(H)$
in which a fixed token $T_1$ is located at $u$ and $v$, respectively.
Along every walk in $H$, the token $T_1$ can only move within its component
of $G$, so ${\boldsymbol{x}}$ and ${\boldsymbol{y}}$ cannot be joined, a contradiction. Hence $G$ is
connected.

Moreover, if $k=ns$, then every vertex of $G$ is full in every configuration,
so no token can move. Since $H$ has more than one vertex, it is disconnected.
Thus $k<ns$.

Conversely, suppose that $G$ is connected and $k<ns$. We show that any two
configurations ${\boldsymbol{x}}$ and ${\boldsymbol{y}}$ are joined in $H$.

We first show that ${\boldsymbol{x}}$ can be transformed into a configuration
${\boldsymbol{z}}$ satisfying
\[
a_{\boldsymbol{z}}(v)=a_{\boldsymbol{y}}(v)\qquad\text{for every }v\in V(G).
\]

We use the following simple observation. Let
$v_0v_1\cdots v_\ell$
be a path in $G$ such that $a_W(v_0)>0$ and $a_W(v_\ell)<s$.
Then one token can be transferred from $v_0$ to $v_\ell$, by legal moves,
without changing the occupancies of the intermediate vertices. Indeed, this
follows by induction on $\ell$. If $\ell=1$, move a token directly from
$v_0$ to $v_1$. For $\ell\ge2$, if $v_1$ is not full, first move a token
from $v_0$ to $v_1$ and apply the induction hypothesis to
$v_1\cdots v_\ell$. If $v_1$ is full, first apply the induction hypothesis
to $v_1\cdots v_\ell$, thereby creating a vacancy at $v_1$, and then move
a token from $v_0$ to $v_1$.

If the occupancy vectors of ${\boldsymbol{x}}$ and ${\boldsymbol{y}}$ differ, choose vertices $u,v$ such
that
\[
a_{\boldsymbol{x}}(u)>a_{\boldsymbol{y}}(u),\qquad a_{\boldsymbol{x}}(v)<a_{\boldsymbol{y}}(v).
\]
Since $G$ is connected, there is a path from $u$ to $v$, and the preceding
observation transfers one token from $u$ to $v$. This decreases
\[
\sum_{w\in V(G)}|a_{\boldsymbol{x}}(w)-a_{\boldsymbol{y}}(w)|
\]
by $2$. Repeating the procedure therefore yields a configuration ${\boldsymbol{z}}$
reachable from ${\boldsymbol{x}}$ with the same occupancy vector as ${\boldsymbol{y}}$.

It remains to rearrange the token labels. We show that the positions of any
two tokens can be interchanged by a walk in $H$. Let $T_a$ and $T_b$ be
two tokens. If they occupy the same vertex, there is nothing to prove.
Otherwise, suppose that they are at distinct vertices $u$ and $v$.

We shall need to create vacancies while keeping $T_a$ and $T_b$ fixed.
If a vertex $w$ is full, then, since $k<ns$, some vertex of $G$ is not full.
Choose a shortest path
$w_0,w_1,\ldots,w_t$
from $w$ to a non-full vertex, where $w_0=w$. Then $w_0,\ldots,w_{t-1}$ are full.
For $j=t-1,t-2,\ldots,0$, move a token from $w_j$ to $w_{j+1}$.
At most one of $T_a,T_b$ lies at $w_j$, and $s\ge2$, so a token different
from $T_a,T_b$ can always be chosen. Thus a vacancy can be created at $w$
without moving either distinguished token.

Now choose a path from $u$ to $v$. Move $T_a$ successively along this path
towards $v$. Whenever its next vertex is full, first create a vacancy there
as above. Eventually $T_a$ and $T_b$ occupy the same vertex. Let
${\boldsymbol{x}}_0,{\boldsymbol{x}}_1,\ldots,{\boldsymbol{x}}_t$
be the resulting walk, so that $T_a$ and $T_b$ coincide in ${\boldsymbol{x}}_t$.

Let $\sigma$ be the permutation of the token labels that interchanges
$T_a$ and $T_b$ and fixes all other tokens. Since the tokens are
distinguishable but obey identical movement and capacity rules, $\sigma$
induces an automorphism of $H$. Moreover,
$\sigma({\boldsymbol{x}}_t)={\boldsymbol{x}}_t$.
Hence,
\[
{\boldsymbol{x}}_0,{\boldsymbol{x}}_1,\ldots,{\boldsymbol{x}}_t,
\sigma({\boldsymbol{x}}_{t-1}),\ldots,\sigma({\boldsymbol{x}}_0)
\]
is a walk from ${\boldsymbol{x}}_0$ to $\sigma({\boldsymbol{x}}_0)$, and therefore interchanges the positions
of $T_a$ and $T_b$ while restoring all other tokens to their original
positions.

Finally, since ${\boldsymbol{z}}$ and ${\boldsymbol{y}}$ have the same number of tokens at every vertex, the assignment of token labels in ${\boldsymbol{z}}$ can be transformed into
that of ${\boldsymbol{y}}$ by a sequence of transpositions. By the preceding argument,
each such transposition can be realized by a walk in $H$. Hence, ${\boldsymbol{z}}$ is
connected to ${\boldsymbol{y}}$ and, therefore, ${\boldsymbol{x}}$ is connected to ${\boldsymbol{y}}$.

Since ${\boldsymbol{x}}$ and ${\boldsymbol{y}}$ were arbitrary, $H$ is connected.
\end{proof}

For a graph $H$ and $W\subseteq V(H)$, we denote by $H[W]$ the subgraph of $H$
induced by $W$.

\begin{theorem}\label{th:t-connected}
Let $G$ be a $t$-vertex-connected graph of order $n$, where $t\ge1$,
and let $s,k$ be integers satisfying
$2\le s<k<ns$.
Then $F^s_{k\times1}(G)$ is $t$-vertex-connected.
\end{theorem}

\begin{proof}
Let
$H=F^s_{k\times1}(G)$,
and fix a token $T_1$. For each $v\in V(G)$, let
\[
V_v=\{{\boldsymbol{x}}\in V(H): T_1\text{ is at }v\},
\qquad
H_v=H[V_v].
\]
Thus, the sets $V_v$, $v\in V(G)$, form a partition of $V(H)$.

If $t=1$, then $G$ is connected and hence $H$ is connected by
Theorem~\ref{th:connected}. Thus, suppose that $t\ge2$.

We first establish two claims.

\medskip
\noindent
\textbf{Claim 1.} Each layer $H_v$ is connected.

\smallskip
Since $t\ge2$, the graph $G-v$ is connected. Put
$q=k-s$,
so that
$1\le q<s(n-1)$.
For every $(s-1)$-subset
$A\subseteq\{T_2,\ldots,T_k\}$, let $H_v(A)$ be the subgraph of $H_v$
induced by the configurations in which the tokens at $v$ are precisely
$\{T_1\}\cup A$.

The remaining $q$ tokens lie on $G-v$. If $q\le s$, the capacity condition
does not restrict their motion, and they can be moved independently on the
connected graph $G-v$. If $q>s$, then
\[
H_v(A)\cong F^s_{q\times1}(G-v),
\]
which is connected by Theorem~\ref{th:connected}. Hence, every $H_v(A)$ is
connected.

We next show that these subgraphs are mutually connected. Suppose
\[
A'=(A\setminus\{a\})\cup\{b\}.
\]
Choose a neighbor $w$ of $v$. There is a configuration in $H_v(A)$ in which
$b$ is at $w$ and $w$ has a vacant place: indeed,
\[
q-1\le s(n-1)-2,
\]
so the remaining tokens can be distributed accordingly. Moving first
$a$ from $v$ to $w$, and then $b$ from $w$ to $v$, gives a path of length
two from $H_v(A)$ to $H_v(A')$. Since any two $(s-1)$-subsets are joined
by a sequence of single-element replacements, all configurations in which
$v$ is full lie in the same component of $H_v$.

Finally, let ${\boldsymbol{x}}\in V_v$ and suppose that $v$ is not full in ${\boldsymbol{x}}$. Since
$k>s$, some token lies outside $v$. Choose an occupied vertex $w\ne v$
at minimum distance from $v$, and let
\[
v=v_0,v_1,\ldots,v_\ell=w
\]
be a shortest path. By the choice of $w$, the internal vertices
$v_1,\ldots,v_{\ell-1}$ are empty. Hence, a token at $w$ can be moved along
this path to $v$, keeping $T_1$ fixed. Repeating this procedure eventually
fills $v$. Thus every vertex of $H_v$ is connected to a configuration in
one of the subgraphs $H_v(A)$ and, therefore, $H_v$ is connected.

\medskip
\noindent
\textbf{Claim 2.} For every $u\in V(G)$ and every ${\boldsymbol{x}}\in V_u$, there exist
$t$ paths in $H$, meeting only at ${\boldsymbol{x}}$, whose internal vertices lie in
$V_u$ and whose terminal vertices lie in $t$ distinct layers
$V_{u_1},\ldots,V_{u_t}$, where $u_1,\ldots,u_t$ are distinct neighbors
of $u$.

\smallskip
First suppose that some vertex $z\ne u$ is not full in ${\boldsymbol{x}}$. By Menger's
theorem, there are $t$ internally vertex-disjoint $u$--$z$ paths in $G$.
Their first vertices after $u$ are distinct; denote them by
$u_1,\ldots,u_t$.

For each such path, truncate it at the first non-full vertex after $u$:
\[
u=w_0, w_1,\ldots,w_\ell.
\]
Thus $w_1,\ldots,w_{\ell-1}$ are full and $w_\ell$ is not full. Moving one
token successively
\[
w_{\ell-1}\to w_\ell,\quad
w_{\ell-2}\to w_{\ell-1},\quad\ldots,\quad
w_1\to w_2
\]
creates a vacancy at $w_1$. We may then move $T_1$ from $u$ to $w_1$.
All preceding moves keep $T_1$ at $u$. Since the original $u$--$z$ paths
are internally disjoint, the $t$ resulting paths in $H$ meet only at ${\boldsymbol{x}}$,
and their terminal vertices belong to the distinct layers $V_{u_i}$.

Now suppose that every vertex other than $u$ is full in ${\boldsymbol{x}}$. Then $u$ is
not full, since $k<ns$. As $\delta(G)\ge t$ (where $\delta$ is the minimum degree) choose distinct neighbors
$u_1,\ldots,u_t$ of $u$. For each $i$, move one token from $u_i$ to $u$,
and then move $T_1$ from $u$ to $u_i$. These $t$ paths have length two,
meet only at ${\boldsymbol{x}}$, and terminate in the distinct layers $V_{u_i}$.
This proves the claim.

\medskip
We now prove the theorem. Let $D\subseteq V(H)$ with
\[
|D|\le t-1.
\]
We show that $H-D$ is connected. Define
\[
U=\{v\in V(G):D\cap V_v\ne\emptyset\}.
\]
Since
$|U|\le |D|<t,$
the graph $G-U$ is connected.

Consider
\[
K=
H\left[
\bigcup_{v\in V(G)\setminus U}V_v
\right].
\]
By Claim~1, each layer $H_v$ with $v\notin U$ is connected. Moreover, if
$uv\in E(G-U)$, then some edge of $H$ joins $V_u$ and $V_v$: place $T_1$
at $u$ and distribute the remaining tokens while leaving a vacant place at
$v$, which is possible because
\[
k-1\le ns-2,
\]
and then move $T_1$ from $u$ to $v$. Since $G-U$ is connected, it follows
that $K$ is connected.

It remains to connect to $K$ the vertices lying in layers indexed by $U$.
Let
\[
{\boldsymbol{x}}\in V(H)\setminus D,\qquad {\boldsymbol{x}}\in V_u,\quad u\in U.
\]
By Claim~2, there are $t$ paths from ${\boldsymbol{x}}$ whose internal vertices lie in
$V_u$, which meet only at ${\boldsymbol{x}}$, and whose terminal vertices lie in $t$
distinct neighboring layers.

At most $|D\cap V_u|$ of these paths contain a deleted internal vertex.
Moreover, at most $|U|-1$ terminate in a layer indexed by
$U\setminus\{u\}$. Since every such layer contains at least one vertex of
$D$,
\[
|D\cap V_u|+|U|-1
\le
|D\cap V_u|
+\sum_{v\in U\setminus\{u\}}|D\cap V_v|
=|D|
\le t-1.
\]
Hence, at least one of the $t$ paths avoids $D$ and terminates in a layer
$V_v$ with $v\notin U$ and, therefore, in $K$. Thus ${\boldsymbol{x}}$ is connected to
$K$ in $H-D$.

Consequently, every vertex of $H-D$ lies in the same component as $K$.
Hence, $H-D$ is connected whenever $|D|\le t-1$ and, therefore, $H$ is
$t$-vertex-connected.
\end{proof}

The inequality can be strict: for example,
$\kappa(K_{1,3})=1$, whereas
$\kappa(F^2_{3\times1}(K_{1,3}))\ge2$.


\section{Metric properties}
\label{sec:distances}

\subsection{General case}
Throughout this section, $G$ is a simple and connected graph, $n=|V(G)|$, and
$2\le s<k\le ns$. The case $s=1$ was studied
in~Dalf\'o, Fern\'andez, Fiol, and Reyes \cite{dffr26}. 
We define the diameter of a disconnected graph to be $+\infty$. 
A \emph{vacant place} is one unit of unused capacity.
For configurations ${\boldsymbol{x}}=(x_1,\ldots,x_k)$ and ${\boldsymbol{y}}=(y_1,\ldots,y_k)$,
put
$$
L_G({\boldsymbol{x}},{\boldsymbol{y}})=\sum_{i=1}^k d_G(x_i,y_i).
$$
Since every labeled token must reach its destination,
$$
d_{F^s_{k\times1}(G)}({\boldsymbol{x}},{\boldsymbol{y}})\ge L_G({\boldsymbol{x}},{\boldsymbol{y}}).
$$
Note that capacity constraints may force additional moves.
\begin{proposition}
Let $G$ be a connected graph with at least two vertices, and let
$2\le s<k<2s$. Then
\[
\diam\bigl(F^s_{k\times1}(G)\bigr)=k\diam(G).
\]
\end{proposition}

\begin{proof}
Write $H=F^s_{k\times1}(G)$ and $d=\diam(G)$. Let
\[
{\boldsymbol{x}}=(x_1,\ldots,x_k),\qquad {\boldsymbol{y}}=(y_1,\ldots,y_k)
\]
be two configurations of $H$, and choose, for each token $T_i$, a shortest
$x_i$--$y_i$ path in $G$. Put
\[
L_G({\boldsymbol{x}},{\boldsymbol{y}})=\sum_{i=1}^k d_G(x_i,y_i).
\]
Since $k<2s$, at most one vertex of $G$ can be full in any configuration.

We first observe that, unless a full vertex contains only tokens already at
their destinations, some unfinished token (that is, a token that has not yet
reached its destination) can be moved one step along its
chosen shortest path. Indeed, if an unfinished token is not full,
the move is possible. Otherwise, its next vertex is the unique full
vertex $v$. Since not all tokens at $v$ are finished, some unfinished token
at $v$ can be moved along its chosen path; its next vertex is different from
$v$ and, hence, is not full.

Thus, starting from ${\boldsymbol{x}}$, we can successively move tokens along their chosen
shortest paths until either all tokens reach their destinations or we are
about to create a full vertex $v$ consisting entirely of finished tokens.
In the latter case, let $T_i$ be the token whose next move would fill $v$.
Then $v$ already contains $s-1$ finished tokens, while $T_i$ is at a
neighbor of $v$ and has destination $v$. Postpone this last move.

Keep these $s-1$ tokens and $T_i$ fixed. There remain $k-s<s$ other tokens.
They can all be moved to their destinations along their chosen shortest
paths. Indeed, if no vertex is full, any unfinished token can move. If $v$ is
full, its additional token is one of these remaining tokens and is
unfinished, so it can leave $v$. If some vertex $w\ne v$ is full, then,
because only
$k-s+1\le s$
tokens lie outside the fixed $s-1$ tokens at $v$, all of them must be at
$w$; hence, unless all remaining tokens are finished, an unfinished one at
$w$ can leave. In every case, one of the remaining shortest routes can be
shortened. After these tokens have reached their destinations, move $T_i$
to $v$. Hence,
\[
d_H({\boldsymbol{x}},{\boldsymbol{y}})\le L_G({\boldsymbol{x}},{\boldsymbol{y}})
\]
unless ${\boldsymbol{x}}$ already contains a full vertex consisting entirely of finished
tokens.

In this exceptional case, let $v$ be that vertex. Temporarily move one of
its $s$ tokens to a neighbor of $v$, apply the preceding argument, and move
it back at the end. Since the $s$ tokens initially at $v$ are already at
their destinations,
\[
L_G({\boldsymbol{x}},{\boldsymbol{y}})\le (k-s)d,
\]
and, therefore,
\[
d_H({\boldsymbol{x}},{\boldsymbol{y}})\le L_G({\boldsymbol{x}},{\boldsymbol{y}})+2
   \le (k-s)d+2
   \le kd,
\]
because $s\ge2$ and $d\ge1$. Consequently,
\[
\diam(H)\le kd.
\]

For the reverse inequality, choose $u,v\in V(G)$ with $d_G(u,v)=d$.
Let ${\boldsymbol{x}}$ have $s$ tokens at $u$ and the remaining $k-s$ tokens at $v$, and
let ${\boldsymbol{y}}$ be obtained by sending each of these tokens to the opposite vertex.
Both configurations are admissible since $k-s<s$. Every token must travel
distance $d$, so
\[
d_H({\boldsymbol{x}},{\boldsymbol{y}})\ge L_G({\boldsymbol{x}},{\boldsymbol{y}})=kd.
\]
Thus, $\diam\bigl(F^s_{k\times1}(G)\bigr)=kd$.
\end{proof}

\subsection{Complete graphs}

\begin{theorem} 
\label{thm:diam-complete-short}
Let $n\ge2$ and $k<ns$. Let $m=|\{i:x_i\ne y_i\}|$, and let $c$
count the components of $D({\boldsymbol{x}},{\boldsymbol{y}})$ whose vertices are all full in ${\boldsymbol{x}}$.
Then
$$
d_{F^s_{k\times1}(K_n)}({\boldsymbol{x}},{\boldsymbol{y}})=m+c.
$$
Moreover,
$$
\operatorname{diam}(F^s_{k\times1}(K_n))=
\begin{cases}
k+\left\lfloor\dfrac{k}{2s}\right\rfloor
 &\text{for }k\le(n-1)s,\\[3mm]
k+\left\lfloor\dfrac n2\right\rfloor-1
 &\text{for }(n-1)s<k<ns.
\end{cases}
$$
\end{theorem}
\begin{proof}
A component full in ${\boldsymbol{x}}$ is also full in ${\boldsymbol{y}}$, since its total number
of tokens is unchanged. Its first move must send a token outside
the component, although that token must eventually return. Thus, 
at least one additional move is needed for each such component,
giving the lower bound $m+c$.
For a component with spare capacity, let $a_v,b_v$ be its token
counts in ${\boldsymbol{x}},{\boldsymbol{y}}$. Reverse its arcs, add a vertex $z$, and add
$s-a_v$ arcs $z\to v$ and $s-b_v$ arcs $v\to z$. 
The resulting directed multigraph is balanced, that is, every vertex has
equal in-degree and out-degree. Indeed, at each $v$ the difference between
out-degree and in-degree is
\[
(b_v-a_v)+(s-b_v)-(s-a_v)=0,
\]
and the same holds at $z$ since $\sum_v a_v=\sum_v b_v$.
Hence, each connected component has an Euler circuit.
Each reversed request arc moves its token into the current vacant
place; arcs incident with $z$ start or end a vacancy sequence.
Thus, every requested token moves once. For a full component,
move one requested token temporarily to an external non-full vertex,
process the remaining arcs of an Euler circuit backward, and
finally deliver the stored token. This uses one extra move and
restores the external vertex, proving the distance formula.
Each full component contains at least two vertices and $2s$ tokens,
and some non-full vertex lies outside them. Hence,
$$
c\le\min\left\{
\left\lfloor\frac{k}{2s}\right\rfloor,
\left\lfloor\frac{n-1}{2}\right\rfloor
\right\}.
$$
This gives both diameter bounds except when $n=2p+1$ and
$(n-1)s<k<ns$. In that case, $c=p$ forces $p$ full pairs and at
least one token with $x_i=y_i$ at the remaining vertex; therefore
$m+c\le k+p-1$ as well.
For equality in the first range, exchange full token groups on
$\lfloor k/(2s)\rfloor$ disjoint pairs of vertices. Distribute any
remaining tokens among the remaining vertices and cyclically shift
their vertex names; at least two vertices remain when needed.
For the second range, use $\lfloor n/2\rfloor-1$ full pairs and
cyclically shift the remaining two vertices when $n$ is even,
or three when $n$ is odd. Distribute the remaining tokens within
these vertices respecting capacity. In both constructions, all
tokens change position, and exactly the required number of full
components occurs.
\end{proof}

\subsection{Stars}

\begin{theorem}
\label{thm:diam-stars-short}
Let $S_n=K_{1,n}$, with $n\ge2$, and $k<(n+1)s$. Then
$$
\operatorname{diam}(F^s_{k\times1}(S_n))=
\begin{cases}
2\min\{k,ns\}&\text{for }k\le(n+1)s-2,\\[2mm]
2ns+2\left\lfloor\dfrac n2\right\rfloor
 &\text{for }k=(n+1)s-1.
\end{cases}
$$
For $S_1=K_2$, the diameter is $k$ when $k<2s$.
\end{theorem}
\begin{proof}
Let $0$ be the center of $S_n$, and let
$a=a_{\boldsymbol{x}}(0)$ and $b=a_{\boldsymbol{y}}(0)$
be the numbers of tokens at the center in two configurations ${\boldsymbol{x}}$ and ${\boldsymbol{y}}$.
Recall that
\[
L_{S_n}({\boldsymbol{x}},{\boldsymbol{y}})=\sum_{i=1}^k d_{S_n}(x_i,y_i).
\]
Since a token starting or ending at the center travels at most one edge,
whereas any other token travels at most two,
\begin{equation}\label{eq:star-L}
L_{S_n}({\boldsymbol{x}},{\boldsymbol{y}})\le 2k-a-b.
\end{equation}

We distinguish two cases.

\medskip
\noindent
\textbf{Case 1:} $k=(n+1)s-1$.
There is now exactly one vacant place. First, suppose that this vacancy is at
the center in both ${\boldsymbol{x}}$ and ${\boldsymbol{y}}$. Thus, the center contains $s-1$ tokens and
every leaf contains $s$ tokens in both configurations.

Let $m_L$ be the number of misplaced tokens that start at leaves, and let
$c_L$ be the number of nontrivial components of $D({\boldsymbol{x}},{\boldsymbol{y}})$ not containing
the center. Since the initial and final occupancies agree at every vertex,
each component of $D({\boldsymbol{x}},{\boldsymbol{y}})$ is balanced.

A pair of moves
\[
\text{leaf}\longrightarrow 0\longrightarrow\text{leaf}
\]
exchanges a leaf token through the vacant central place and restores the
vacancy to the center. Every misplaced leaf token requires one such
exchange. Moreover, each component not containing the center requires one
additional exchange in order to introduce a temporary central token and
restore it afterward. Conversely, following an Euler circuit in each
component realizes all requests with exactly these exchanges. Hence, 
\[
d_{F^s_{k\times1}(S_n)}({\boldsymbol{x}},{\boldsymbol{y}})=2(m_L+c_L). 
\]

Clearly,
\[
m_L\le ns,
\qquad
c_L\le \left\lfloor\frac n2\right\rfloor,
\]
since every component counted by $c_L$ contains at least two leaves.
Therefore,
\[
d_{F^s_{k\times1}(S_n)}({\boldsymbol{x}},{\boldsymbol{y}})
\le
2ns+2\left\lfloor\frac n2\right\rfloor
=:D_0.
\]

Equality is attained by keeping the central tokens fixed and exchanging the
$s$-token sets of the leaves in pairs; when $n$ is odd, use one cyclically
permuted triple of leaves and pair the remaining leaves. Then
\[
m_L=ns,\qquad
c_L=\left\lfloor\frac n2\right\rfloor,
\]
so (5) gives distance $D_0$.

It remains to consider the cases not covered above. We shall repeatedly use
the following observation. If one token is moved temporarily from the center
to a non-full leaf, producing a configuration ${\boldsymbol{x}}'$, then
\[
d({\boldsymbol{x}},{\boldsymbol{x}}')=1
\qquad\text{and}\qquad
L_{S_n}({\boldsymbol{x}}',{\boldsymbol{y}})\le L_{S_n}({\boldsymbol{x}},{\boldsymbol{y}})+1.
\]
Thus, if the pair $({\boldsymbol{x}}',{\boldsymbol{y}})$ can be realized using only shortest routes, then
\begin{equation}\label{eq:star-temp}
d({\boldsymbol{x}},{\boldsymbol{y}})\le L_{S_n}({\boldsymbol{x}},{\boldsymbol{y}})+2.
\end{equation}

Suppose first that $(a,b)=(s-1,s-1)$. Since there are at least two vacant
places altogether, some leaf is non-full. Move one central token temporarily
to such a leaf. The new central occupancy is $s-2$, so the preceding case
applies and~\eqref{eq:star-temp} follows.

The same argument applies if $a=s$ and $b\le s-2$: move one central token
to a non-full leaf. The new pair of central occupancies is
$(s-1,b)$, which is covered by the preceding case. Hence again
$d({\boldsymbol{x}},{\boldsymbol{y}})\le L_{S_n}({\boldsymbol{x}},{\boldsymbol{y}})+2$.
By symmetry, the same conclusion holds when $b=s$ and $a\le s-2$.

Now suppose that $(a,b)=(s,s-1)$. Since the center has one more token in
${\boldsymbol{x}}$ than in ${\boldsymbol{y}}$, the request digraph $D({\boldsymbol{x}},{\boldsymbol{y}})$ contains a directed request
path starting at the center and ending at a deficit leaf. As ${\boldsymbol{x}}$ has at
least two vacant places, all of them on leaves, we may move a central token
different from the token initiating this path to a leaf while leaving a
vacant place at the terminal leaf of the path. Process the request path
backward, using shortest routes. After doing so, the center contains
$s-2$ tokens, and the remaining requests fall under the preceding case.
The preliminary displacement contributes at most two moves beyond
$L_{S_n}({\boldsymbol{x}},{\boldsymbol{y}})$ and, hence,
$d({\boldsymbol{x}},{\boldsymbol{y}})\le L_{S_n}({\boldsymbol{x}},{\boldsymbol{y}})+2$.

The case $(a,b)=(s-1,s)$ follows by reversing the roles of ${\boldsymbol{x}}$ and ${\boldsymbol{y}}$.

Finally, suppose that $(a,b)=(s,s)$. Move one central token of ${\boldsymbol{x}}$
temporarily to a non-full leaf, obtaining ${\boldsymbol{x}}'$ with central occupancy
$s-1$. By the preceding case,
$d({\boldsymbol{x}}',{\boldsymbol{y}})\le L_{S_n}({\boldsymbol{x}}',{\boldsymbol{y}})+2$.
Since $d({\boldsymbol{x}},{\boldsymbol{x}}')=1$ and
$L_{S_n}({\boldsymbol{x}}',{\boldsymbol{y}})\le L_{S_n}({\boldsymbol{x}},{\boldsymbol{y}})+1$, we obtain
\[
d({\boldsymbol{x}},{\boldsymbol{y}})\le L_{S_n}({\boldsymbol{x}},{\boldsymbol{y}})+4.
\]

Thus, the extra cost beyond $L_{S_n}({\boldsymbol{x}},{\boldsymbol{y}})$ is zero in the basic case,
at most $2$ in all the remaining cases except $(a,b)=(s,s)$, and at most
$4$ in that case. Since $a,b\ge r$ and $r\le s-2$, these bounds are in
every case at most
$$
a+b-2r.
$$
Consequently, by~\eqref{eq:star-L},
\[
\begin{aligned}
d_{F^s_{k\times1}(S_n)}({\boldsymbol{x}},{\boldsymbol{y}})
&\le L_{S_n}({\boldsymbol{x}},{\boldsymbol{y}})+a+b-2r\\
&\le 2k-2r\\
&=2\min\{k,ns\}.
\end{aligned}
\]

\textbf{Case 2:} $k\le (n+1)s-2$.
Put
\[
r=\max\{0,k-ns\}.
\]
Since the leaves have total capacity $ns$, every configuration has at least
$r$ tokens at the center. Moreover, the assumption on $k$ gives
$r\le s-2$.

We first consider the case
\[
a,b\le s-1,\qquad \min\{a,b\}\le s-2.
\]
Let $D({\boldsymbol{x}},{\boldsymbol{y}})$ be the directed multigraph having one arc $x_i\to y_i$
for each token with $x_i\ne y_i$. Its arcs may be decomposed into directed
paths from surplus to deficit vertices and directed cycles. If $a\ne b$,
process first the paths incident with the center, starting from the endpoint
having the larger central occupancy. This may be done along shortest routes
until both remaining configurations have $\min\{a,b\}\le s-2$ tokens at the
center.

There are then at least two vacant places at the center. Process the
remaining directed paths in reverse order. A directed cycle containing only
leaves is handled by temporarily holding one token at the center and then
processing the remaining requests around the cycle; a cycle containing the
center is treated similarly. Since two central places are available, every
token can follow a shortest route. Hence, in this case,
\[
d_{F^s_{k\times1}(S_n)}({\boldsymbol{x}},{\boldsymbol{y}})=L_{S_n}({\boldsymbol{x}},{\boldsymbol{y}}).
\]

The remaining possibilities occur when one or both center occupancies are
$s-1$ or $s$. They reduce to the preceding case by temporarily moving a
central token to a non-full leaf. Such a temporary displacement and its
subsequent correction require at most two additional moves. Thus the extra
cost beyond $L_{S_n}({\boldsymbol{x}},{\boldsymbol{y}})$ is at most $2$, except when $a=b=s$, when at
most two such displacements are required and the extra cost is at most $4$.
Consequently, in every case the additional cost is at most
\[
a+b-2r,
\]
because $a,b\ge r$ and $r\le s-2$. Therefore, by (4),
$$
d_{F^s_{k\times1}(S_n)}(X,{\boldsymbol{y}})
\le L_{S_n}({\boldsymbol{x}},{\boldsymbol{y}})+a+b-2r
\le 2k-2r
 =2\min\{k,ns\}.
$$

To attain this bound, if $k\le ns$, place all $k$ tokens on the leaves and
cyclically shift the leaf containing each token. If $k>ns$, put
$r=k-ns$ tokens at the center and keep them fixed, fill all leaves, and
cyclically shift the $ns$ leaf tokens. In either case,
$\min\{k,ns\}$ tokens travel between distinct leaves and, hence, each has
distance $2$. Thus
\[
\diam\bigl(F^s_{k\times1}(S_n)\bigr)
   =2\min\{k,ns\}.
\]





\end{proof}

\subsection{Odd cycles}

Let $n=2m+1\ge3$. We first establish a general lower bound for the diameter over odd cycles.

\begin{proposition}
\label{prop:odd-cycle-diameter-lower}
Let $n=2m+1\ge3$ and $2\le s<k<ns$. Then
$$
\operatorname{diam}(F^s_{k\times1}(C_{2m+1}))\ge
\begin{cases}
km+\left\lfloor\dfrac{k}{2s}\right\rfloor
 &\text{for }k\le2ms,\\[2mm]
km&\text{for }2ms<k<(2m+1)s.
\end{cases}
$$
\end{proposition}
\begin{proof}
Label the cycle with vertices $0,\ldots,2m$ modulo $2m+1$.
Rotating every token by $m$ vertices gives the lower bound $km$.
For $k\le2ms$, write
$$
k=2rs+t,\qquad r=\left\lfloor\frac{k}{2s}\right\rfloor,
\qquad \text{for }0\le t<2s.
$$
Assume $r\ge1$, since $r=0$ is already covered. Fill the vertices
$I=\{0,\ldots,2r-1\}$. If $t>0$, then $r<m$; place the remaining
tokens at $2r,2r+1$, at most $s$ at each. Send tokens initially
at $j<r$ to $j+m$, and all others to $j-m$.
Two initially occupied vertices from the two groups cannot be assigned
the same destination vertex. Indeed, such an equality would require
$$
i+m\equiv j-m\mod(2m+1),
$$
that is,
$$
j-i\equiv2m\mod(2m+1),
$$
which is impossible because $1\le j-i\le2m-1$.
Hence, the final configuration is valid and
$d_{C_{2m+1}}(x_i,y_i)=m$ for every token.
In any legal sequence, mark the first departure from each $j\in I$
by an arrow to the receiving neighbor. If that neighbor lies in
$I$, its first departure occurred earlier, since it was initially
full. Following the arrows from $r-1$, therefore, exits $I$ through
$0$ or $2r-1$, with strictly decreasing departure times.
In the first case, the first departing tokens from $0,\ldots,r-1$
move opposite to their unique shortest routes; in the second,
those from $r,\ldots,2r-1$ do so. These are distinct tokens making
their first moves. Each needs at least $m+1$ moves, since its
distance to its destination after that first move is still $m$.
Thus every sequence has length at least
$$
(k-r)m+r(m+1)=km+r.
$$
\end{proof}
For $C_3=K_3$, Theorem~\ref{thm:diam-complete-short} gives
$$
\operatorname{diam}(F^s_{k\times1}(C_3))=
\begin{cases}
k&\text{for }s<k<2s,\\
k+1&\text{for }k=2s,\\
k&\text{for }2s<k<3s.
\end{cases}
$$
The following formula is, therefore, proved as a lower bound, with
equality established when $s<k<2s$ or $n=3$.
\begin{conjecture}
\label{conj:odd-cycle-diam}
Let $n=2m+1\ge3$ and $2\le s<k<(2m+1)s$. Then
$$
\operatorname{diam}(F^s_{k\times1}(C_{2m+1}))=
\begin{cases}
km+\left\lfloor\dfrac{k}{2s}\right\rfloor,
 &s<k\le2ms,\\[2mm]
km,&2ms<k<(2m+1)s.
\end{cases}
$$
\end{conjecture}

\begin{example}
For \(F^2_{4\times1}(C_5)\), the pair
\({\boldsymbol{x}}=(1,1,2,2)\), \({\boldsymbol{y}}=(3,3,5,5)\) has \(L_{C_5}({\boldsymbol{x}},{\boldsymbol{y}})=8\), but every realization requires one additional move, so \(d({\boldsymbol{x}},{\boldsymbol{y}})=9\).
\end{example}
The exact computations in Table~\ref{tab:odd-cycle-comput} agree
with Conjecture~\ref{conj:odd-cycle-diam} in every tested case.
The corresponding pseudo-code is given in the Appendix.

\begin{table}[ht]
\centering
\caption{Exact computational values of
$\operatorname{diam}(F^s_{k\times1}(C_n))$ for odd cycles.}
\label{tab:odd-cycle-comput}
\renewcommand{\arraystretch}{1.15}
\begin{tabular}{|c|c|c|c|}
\hline
$n$ & $s$ & Values of $k$ &
$\operatorname{diam}(F^s_{k\times1}(C_n))$ \\
\hline\hline
$5$ & $2$ &
$3,4,5,6,7,8,9$ &
$6,9,11,13,15,18,18$ \\
$5$ & $3$ &
$4,5,6,7,8,9,10,11,12,13,14$ &
$8,10,13,15,17,19,21,23,26,26,28$ \\
$5$ & $4$ &
$8,12,16,17$ &
$17,25,34,34$ \\
$7$ & $2$ &
$3,4,5,6,7,8,9,10,11,12,13$ &
$9,13,16,19,22,26,29,32,35,39,39$ \\
$7$ & $3$ &
$6,7,8$ &
$19,22,25$ \\
$9$ & $2$ &
$3,4,5,6$ &
$12,17,21,25$ \\
\hline
\end{tabular}
\end{table}



\section{Clique number}
\label{sec:clique}
Throughout this section, $G$ is a finite simple undirected graph. 
\begin{lemma}\label{clique}
Let $G$ be a graph of order $n\ge2$, and let $s, k$ be integers satisfying $1\le s< k\le ns$. 
For every nonempty clique $\mathcal K$ of $F^s_{k\times1}(G)$, there is a coordinate $i$ such that all configurations in $\mathcal K$ agree in every coordinate other than $i$. Moreover, the values occurring in coordinate $i$ form a clique of $G$ of order $|\mathcal K|$.
\end{lemma} 
\begin{proof}
The assertion is immediate if $|\mathcal K|=1$. Suppose that $|\mathcal K|\ge2$, and choose distinct configurations
\[
{\boldsymbol{x}}=(x_1,\ldots,x_k),\qquad
{\boldsymbol{y}}=(y_1,\ldots,y_k)
\]
in $\mathcal K$. Since ${\boldsymbol{x}}$ and ${\boldsymbol{y}}$ are adjacent, there is a unique index $i\in\{1,\ldots,k\}$ such that $x_i\ne y_i$, whereas $x_\ell=y_\ell$ for every $\ell\ne i$.
Consider any $Z=(z_1,\ldots,z_k)\in\mathcal K\setminus\{{\boldsymbol{x}},{\boldsymbol{y}}\}$. Since ${\boldsymbol{z}}$ is adjacent to ${\boldsymbol{x}}$, it differs from ${\boldsymbol{x}}$ in exactly one coordinate, say $j$. If $j\ne i$, then
\[
z_i=x_i\ne y_i
\qquad\text{and}\qquad
z_j\ne x_j=y_j.
\]
Thus ${\boldsymbol{z}}$ and ${\boldsymbol{y}}$ differ in at least two coordinates, contradicting their adjacency. Hence, $j=i$.
It follows that all configurations in $\mathcal K$ agree outside coordinate $i$. 
Set
\[
C=\{z_i:Z\in\mathcal K\}.
\]
Since all configurations in $\mathcal K$ agree outside coordinate $i$,
distinct configurations have distinct $i$-th entries. Hence,
$|C|=|\mathcal K|$.
If $c,c'\in C$ are distinct, the corresponding configurations
in $\mathcal K$ differ only in coordinate $i$. Their adjacency
therefore implies that $cc'\in E(G)$. Thus $C$ is a clique of $G$.
\end{proof}
\begin{theorem}
Let $G$ be a graph of order $n\ge2$, and let $s, k$ be integers satisfying $1\le s< k\le ns$. Then
\[
\omega\left(F^s_{k\times1}(G)\right)
	=\min\{\omega(G),\,ns-k+1\}.
\]
\end{theorem}
\begin{proof}
Write $H=F^s_{k\times1}(G)$, and label the tokens $T_1,\ldots,T_k$. 
Let $\mathcal K$ be a nonempty clique of $H$ of order $q$. 
By Lemma~\ref{clique}, there is an index $i$ such that the
positions of $T_i$ form a clique $C$ of $G$ of order $q$,
while each of the other $k-1$ tokens occupies the same vertex
in every configuration of $\mathcal K$. 
Hence,
\[
q\le \omega(G).
\]
For every $c\in C$, some configuration in $\mathcal K$ places
$T_i$ at $c$. Therefore, at most $s-1$ fixed tokens are located
at $c$. Each vertex outside $C$ contains at most $s$ fixed tokens.
Consequently,
\[
k-1\le q(s-1)+(n-q)s=ns-q,
\]
and, hence,
\[
q\le ns-k+1.
\]
Thus,
\[
\omega(H)\le \min\{\omega(G),\,ns-k+1\}.
\]
For the reverse inequality, let
\[
r=\min\{\omega(G),\,ns-k+1\},
\]
and choose a clique $C$ of $G$ of order $r$. Since
\[
r(s-1)+(n-r)s=ns-r\ge k-1,
\]
we can place the labeled tokens $T_1,\ldots,T_{k-1}$ so that
each vertex of $C$ receives at most $s-1$ tokens and every vertex
outside $C$ receives at most $s$ tokens. Let $c_i$ be the position
of $T_i$.
For each $v\in C$, define
\[
{\boldsymbol{x}}_v=(c_1,\ldots,c_{k-1},v).
\]
These are valid configurations. If $u,v\in C$ are distinct, then
${\boldsymbol{x}}_u$ and ${\boldsymbol{x}}_v$ differ only in the position of $T_k$, and
$uv\in E(G)$. Hence, ${\boldsymbol{x}}_u{\boldsymbol{x}}_v\in E(H)$. Therefore,
$\{{\boldsymbol{x}}_v:v\in C\}$ is a clique of $H$ of order $r$, so
\[
\omega(H)\ge r.
\]
Combining the two inequalities gives
\[
\omega(H)=\min\{\omega(G),\,ns-k+1\}.
\]
\end{proof}


\section{Chromatic number}
\label{sec:chromatic}
\begin{theorem}\label{chromatic}
Let $G$ be a finite simple undirected graph of order $n\ge2$, and let $s$ and $k$ be integers satisfying $1\le s<k\le ns$. Then
\[
\chi\left(F^s_{k\times1}(G)\right)\le\chi(G).
\]
Moreover, equality holds if there exists a nonempty set $U\subseteq V(G)$ such that
\[
\chi(G[U])=\chi(G)\qquad\text{and}\qquad|U|\le ns-k+1,
\]
where $G[U]$ is the induced subgraph of $G$.
\end{theorem}
\begin{proof}
Write $H=F^s_{k\times1}(G)$, and label the tokens $T_1,\ldots,T_k$. 
If $\chi(G)=1$, then both $G$ and $H$ are edgeless, so $\chi\left(F^s_{k\times1}(G)\right)=\chi(G)=1$. 
Suppose next that $\chi(G)\ge2$. 
Let $c:V(G)\to\{0,\ldots,\chi(G)-1\}$ be a proper coloring of $G$.
For each configuration ${\boldsymbol{x}}=(x_1,\ldots,x_k)\in V(H)$, define
\[
\varphi({\boldsymbol{x}})=\left(\sum_{j=1}^k c(x_j)\right)(\textrm{mod }\chi(G)).
\]
We show that $\varphi$ is a proper coloring of $H$.
Suppose that ${\boldsymbol{x}}=(x_1,\ldots,x_k)$ and ${\boldsymbol{y}}=(y_1,\ldots,y_k)$ are adjacent in $H$.
Then there is a unique index $i$ such that $x_i y_i\in E(G)$, while $x_j=y_j$ for every $j\ne i$. 
Consequently,
\[
\varphi({\boldsymbol{y}})-\varphi({\boldsymbol{x}})
\equiv c(y_i)-c(x_i)\ (\textrm{mod }\chi(G)).
\]
Since $c$ is proper, $c(x_i)\ne c(y_i)$. Both values belong to $\{0,\ldots,\chi(G)-1\}$, so their difference is nonzero modulo $\chi(G)$.
Thus $\varphi({\boldsymbol{x}})\ne\varphi({\boldsymbol{y}})$, proving that
\[
\chi(H)\le \chi(G).
\]
Now suppose that $\emptyset\ne U\subseteq V(G)$ with $|U|\le ns-k+1$ and $\chi(G[U])=\chi(G)$. 
Since
\[
|U|(s-1)+(n-|U|)s=ns-|U|\ge k-1,
\]
we can place $T_1,\ldots,T_{k-1}$ such that each vertex of $U$ contains at most $s-1$ tokens and every vertex outside $U$ contains at most $s$ tokens.
Fix these positions, and let $a_j$ be the position of $T_j$ for $1\le j\le k-1$.
For each $u\in U$, placing $T_k$ at $u$ gives a vertex
\[
{\boldsymbol{x}}_u=(a_1,\ldots,a_{k-1},u)
\]
of $H$.
For distinct $u,v\in U$, we have
\[
{\boldsymbol{x}}_u{\boldsymbol{x}}_v\in E(H)
\quad\Longleftrightarrow\quad
uv\in E(G).
\]
Hence, the subgraph of $H$ induced by $\{{\boldsymbol{x}}_u:u\in U\}$ is isomorphic to $G[U]$.
It follows that
\[
\chi(H)\ge\chi(G[U])=\chi(G).
\]
Together with the upper bound, this gives
$\chi(H)=\chi(G)$.
\end{proof}
\begin{corollary}
Under the assumptions of Theorem~\ref{chromatic}, $\chi\left(F^s_{k\times1}(G)\right)=\chi(G)$ 
if either of the following conditions holds:
\begin{enumerate}
	\item [(i)] $\chi(G)=\omega(G)$ and $k\le ns-\chi(G)+1$;
	\item [(ii)] $s\ge2$ and $k\le n(s-1)+1$. 
\end{enumerate}
\end{corollary}
\begin{proof}
For $(i)$, take $U$ to be the vertex set of a maximum clique of $G$. Then
\[
\chi(G[U])=|U|=\omega(G)=\chi(G), 
\]
and 
\[
|U|=\chi(G)\le ns-k+1.
\]
For $(ii)$, take $U=V(G)$, then 
\[
\chi(G[U])=\chi(G). 
\] 
Note that $k\le n(s-1)+1$ is equivalent to $n\le ns-k+1$. 
Both conclusions follow from Theorem~\ref{chromatic}.
\end{proof}


\section{Planarity over complete base graphs}
\label{sec:planar}

We use the standard planar girth bound
\[
m\le\frac{g}{g-2}(n-2)
\]
and Kuratowski's theorem; see Balakrishnan and Ranganathan \cite{BalakrishnanRanganathan}(Chapter 8).

Here, a
\emph{subdivision} is obtained by replacing edges by internally
 vertex-disjoint paths. 

In this section, we characterize the planarity of supertoken graphs whose base is a complete graph.

\subsection{Planarity over \(K_2\)}
\label{subsec:planarity-K2}

Label the vertices of $K_2$ by $0$ and $1$. Encoding a configuration of $F^s_{k\times1}(K_2)$ by the binary word $\boldsymbol{x}=x_1\cdots x_k$, where $x_i$ is the position of token $t_i$, identifies $F^s_{k\times1}(K_2)$ with an induced subgraph of the hypercube $Q_k$. Indeed, since two vertices of $Q_k$ are adjacent when they differ in exactly one coordinate, $F^s_{k\times1}(K_2)$ is obtained from $Q_k$ by retaining precisely those binary words in which neither $0$ nor $1$ occurs more than $s$ times.

This representation immediately yields the order and size:
$$
\left|V\bigl(F^s_{k\times1}(K_2)\bigr)\right|
=
\begin{cases}
	2^k,
	& 1\leq k\leq s,\\[1mm]
	\displaystyle\sum_{j=k-s}^{s}\binom{k}{j},
	& s<k\leq2s,\\[3mm]
	0,
	& k>2s,
\end{cases} 
$$
$$
\qquad \left|E\bigl(F^s_{k\times1}(K_2)\bigr)\right|
=
\begin{cases}
	k2^{k-1},
	& 1\leq k\leq s,\\[1mm]
	\displaystyle\sum_{j=k-s}^{s-1}(k-j)\binom{k}{j},
	& s<k<2s,\\[3mm]
	0,
	& k\geq2s.
\end{cases}$$
Recall that the \emph{girth} is the length of a shortest cycle. For $s\geq2$ and $2\leq k<2s$, the graph $F^s_{k\times1}(K_2)$ has girth $4$ when $k\leq2s-2$, and girth $6$ when $k=2s-1$.

We now characterize the planarity of $F^s_{k\times1}(K_2)$ using the above formulas for its order, size, and girth together with the planar girth bound.

    \begin{proposition}
\label{prop:K2-girth-edge-bound}
Let $s\geq2$ and $2\leq k<2s$, and let $N$, $M$, and $g$
denote, respectively, the order, size, and girth of
$F^s_{k\times1}(K_2)$. Then
\[
M\leq \frac{g}{g-2}(N-2)
\quad\text{if and only if}\quad
k\in\{2,3\}\ \text{or}\ (s,k)=(3,4).
\]
Moreover, equality holds in all these cases.
\end{proposition}

\begin{proof}
Recall the formulas for $N$, $M$, and $g$ obtained above. We consider two cases.

\medskip
\noindent
\textbf{Case 1: $2\leq k\leq2s-2$.}
In this case, $g=4$, so the required inequality is $M\leq2N-4$. Suppose first that $k\leq s$. Then every binary word of length $k$ is admissible, and hence $N=2^k$, $M=k2^{k-1}$. Thus, the inequality is equivalent to $(4-k)2^{k-1}\geq4$, which holds precisely for $k\in\{2,3\}$, with equality in both cases. Now suppose that $s<k\leq2s-2$, and put $r:=k-s$. Using
$\binom{k}{r}=\binom{k}{s}$ and $(k-j)\binom{k}{j}=(j+1)\binom{k}{j+1}$, the counting formulas give
\[
2M=kN-2r\binom{k}{s},
\qquad
N\geq2\binom{k}{s}.
\]
Consequently, $2M\geq(k-r)N=sN$. If $s\geq4$, then $M\geq2N>2N-4$, so the required inequality fails. It remains to consider $s\leq3$. The restrictions $s<k\leq2s-2$ force $(s,k)=(3,4)$. In this case, $N=14$, $M=24=2N-4$. Thus, this range contributes precisely the pair $(s,k)=(3,4)$.

\medskip
\noindent
\textbf{Case 2: $k=2s-1$.}
In this case, $g=6$, and the only admissible Hamming weights are
$s-1$ and $s$. The counting formulas give $2M=sN$. The required inequality becomes $M\leq\frac32(N-2)$, or equivalently, $(3-s)N\geq6$. This is impossible when $s\geq3$. If $s=2$, then $k=3$ and $N=M=6$, so $M=6=\frac32(N-2)$. 

Combining the two cases proves that the inequality holds precisely
when $k\in\{2,3\}$ or $(s,k)=(3,4)$, and equality holds in each case.
\end{proof}

Proposition~\ref{prop:K2-girth-edge-bound} reduces the planarity
problem to three small parameter families. The completely saturated
case \(k=2s\) must be treated separately, since the resulting graph
has no edges.

\begin{theorem}
	\label{thm:planarity-supertoken-K2}
	Let $s\geq2$ and let $2\leq k\leq2s$. Then $F^s_{k\times1}(K_2)$ 	is planar if and only if
	$
	k\in\{2,3,2s\},
	$
	or
	$
	(s,k)=(3,4).
	$
\end{theorem}

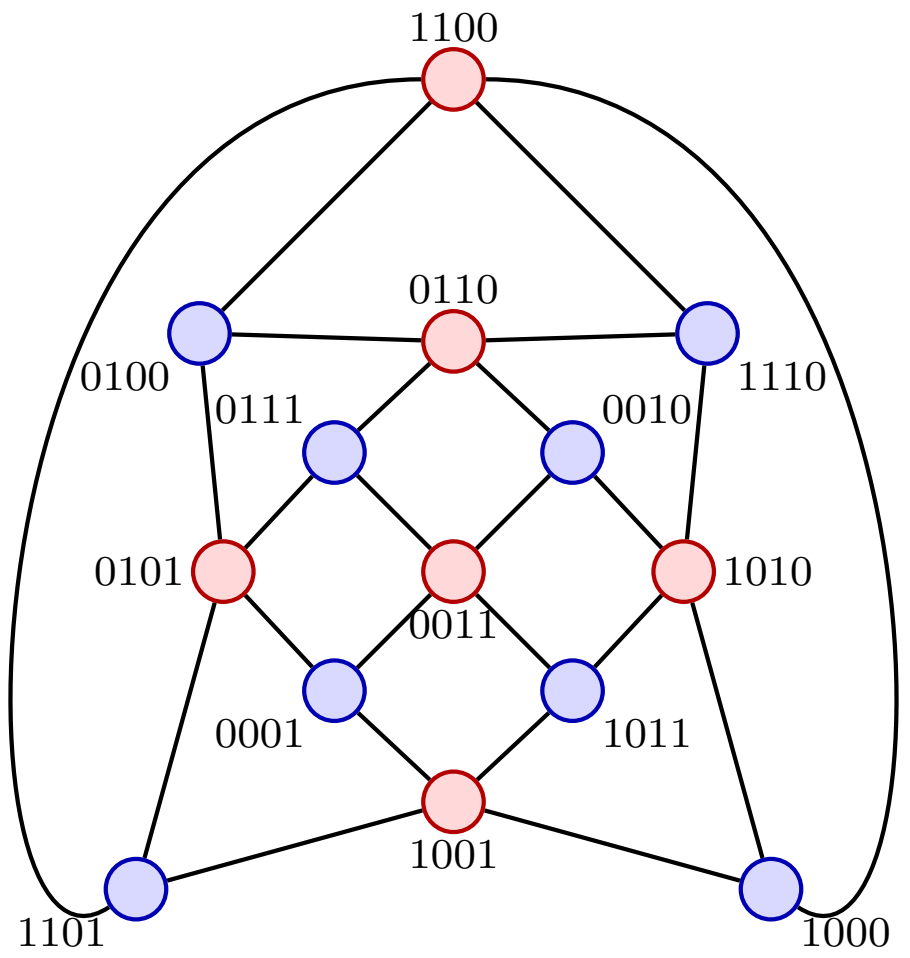
\begin{figure}[ht]
		\centering
		\begin{tikzpicture}[
			scale=1.05,
			oddVertex/.style={
				circle,
				draw=blue!70!black,
				fill=blue!15,
				thick,
				minimum size=4mm,
				inner sep=0pt
			},
			evenVertex/.style={
				circle,
				draw=red!70!black,
				fill=red!15,
				thick,
				minimum size=4mm,
				inner sep=0pt
			},
			vertexLabel/.style={
				font=\scriptsize,
				fill=white,
				inner sep=1pt
			},
			graphEdge/.style={
				black,
				thick
			},
			outerEdge/.style={
				black,
				thick
			}
			]
			
			
			\node[
			oddVertex,
			label={[vertexLabel]below left:$0100$}
			] (o1) at (-1.6,1.5) {};
			
			\node[
			oddVertex,
			label={[vertexLabel]below right:$1110$}
			] (o2) at (1.6,1.5) {};
			
			\node[
			oddVertex,
			label={[vertexLabel]below right:$1000$}
			] (o3) at (2,-2) {};
			
			\node[
			oddVertex,
			label={[vertexLabel]below left:$1101$}
			] (o4) at (-2,-2) {};
			
			\node[
			oddVertex,
			label={[vertexLabel]above left:$0111$}
			] (i1) at (-0.75,0.75) {};
			
			\node[
			oddVertex,
			label={[vertexLabel]above right:$0010$}
			] (i2) at (0.75,0.75) {};
			
			\node[
			oddVertex,
			label={[vertexLabel]below right:$1011$}
			] (i3) at (0.75,-0.75) {};
			
			\node[
			oddVertex,
			label={[vertexLabel]below left:$0001$}
			] (i4) at (-0.75,-0.75) {};
			
			
			\node[
			evenVertex,
			label={[vertexLabel]above:$0110$}
			] (fTop) at (0,1.45) {};
			
			\node[
			evenVertex,
			label={[vertexLabel]right:$1010$}
			] (fRight) at (1.45,0) {};
			
			\node[
			evenVertex,
			label={[vertexLabel]below:$1001$}
			] (fBottom) at (0,-1.45) {};
			
			\node[
			evenVertex,
			label={[vertexLabel]left:$0101$}
			] (fLeft) at (-1.45,0) {};
			
			\node[
			evenVertex,
			label={[vertexLabel]below:$0011$}
			] (fInner) at (0,0) {};
			
			\node[
			evenVertex,
			label={[vertexLabel]above:$1100$}
			] (fOuter) at (0,3.1) {};
			
			
			\draw[graphEdge] (fTop) -- (o1);
			\draw[graphEdge] (fTop) -- (o2);
			\draw[graphEdge] (fTop) -- (i1);
			\draw[graphEdge] (fTop) -- (i2);
			
			
			\draw[graphEdge] (fRight) -- (o2);
			\draw[graphEdge] (fRight) -- (o3);
			\draw[graphEdge] (fRight) -- (i2);
			\draw[graphEdge] (fRight) -- (i3);
			
			
			\draw[graphEdge] (fBottom) -- (o3);
			\draw[graphEdge] (fBottom) -- (o4);
			\draw[graphEdge] (fBottom) -- (i3);
			\draw[graphEdge] (fBottom) -- (i4);
			
			
			\draw[graphEdge] (fLeft) -- (o4);
			\draw[graphEdge] (fLeft) -- (o1);
			\draw[graphEdge] (fLeft) -- (i4);
			\draw[graphEdge] (fLeft) -- (i1);
			
			
			\draw[graphEdge] (fInner) -- (i1);
			\draw[graphEdge] (fInner) -- (i2);
			\draw[graphEdge] (fInner) -- (i3);
			\draw[graphEdge] (fInner) -- (i4);
			
			
			\draw[outerEdge] (fOuter) -- (o1);
			\draw[outerEdge] (fOuter) -- (o2);
			
			\draw[outerEdge]
			(fOuter)
			.. controls (3.2,3.1) and (3.2,-2.8) ..
			(o3);
			
			\draw[outerEdge]
			(fOuter)
			.. controls (-3.2,3.1) and (-3.2,-2.8) ..
			(o4);
			
		\end{tikzpicture}
		
		\caption{A plane embedding of
$F^3_{4\times1}(K_2)\cong Q_4-\{0000,1111\}$, with its bipartition
shown in blue and red.}
		\label{fig:plane-embedding-F3-4-K2}
	\end{figure}

\begin{proof}
Let \(s\geq2\) and \(2\leq k\leq2s\). If $k=2s$, then $F^s_{2s\times1}(K_2)$ is edgeless, and hence planar. Now assume that \(2\leq k<2s\). Let \(N\), \(M\), and \(g\) denote,
respectively, the order, size, and girth of
\(F^s_{k\times1}(K_2)\). If this graph is planar, then the planar
girth bound gives
$
M\leq\frac{g}{g-2}(N-2).
$
By Proposition~\ref{prop:K2-girth-edge-bound}, this is possible only
when $k\in\{2,3\}$ or $(s,k)=(3,4).$ Conversely, all the listed cases are planar. Indeed,
$F^s_{2\times1}(K_2)\cong Q_2\cong C_4$. For $k=3$, we have
$F^s_{3\times1}(K_2)\cong Q_3$ if $s\geq3$, while
$F^2_{3\times1}(K_2)\cong Q_3-\{000,111\}\cong C_6$.
Finally,
$F^3_{4\times1}(K_2)\cong Q_4-\{0000,1111\}$, which is planar by the
embedding in Figure~\ref{fig:plane-embedding-F3-4-K2}.
\end{proof}

\subsection{Planarity for Complete Base Graphs of Order at Least Three}
\label{subsec:planarity-complete-graphs}

We prove that, for \(n\geq3\) and
\(s\geq2\), the graph \(F^s_{k\times1}(K_n)\) is planar only in the
saturated case \(k=ns\), when no token can move.

We first establish several small nonplanar cases that will serve as local obstructions.

\begin{lemma}
	\label{lem:small-nonplanar-supertoken}
	For every $q\in\{2,3,4,5\}$, the graph $F^2_{q\times1}(K_3)$ is nonplanar.
\end{lemma}

\begin{proof}
	Identify $V(K_3)$ with $\{1,2,3\}$ and represent each configuration of
$F^2_{q\times1}(K_3)$ by a word $a_1\cdots a_q\in\{1,2,3\}^q$ in
which each symbol occurs at most twice. Two such words are adjacent
precisely when they differ in one coordinate. 
For each relevant value of $q$, the graph contains a subdivision of $K_{3,3}$,  and is therefore nonplanar by Kuratowski’s theorem. For $q=2$, such a subdivision is displayed in
Figure~\ref{fig:F2-2-K3-K33-subdivision}; the remaining cases are
given in the Appendix.    \end{proof}
	
	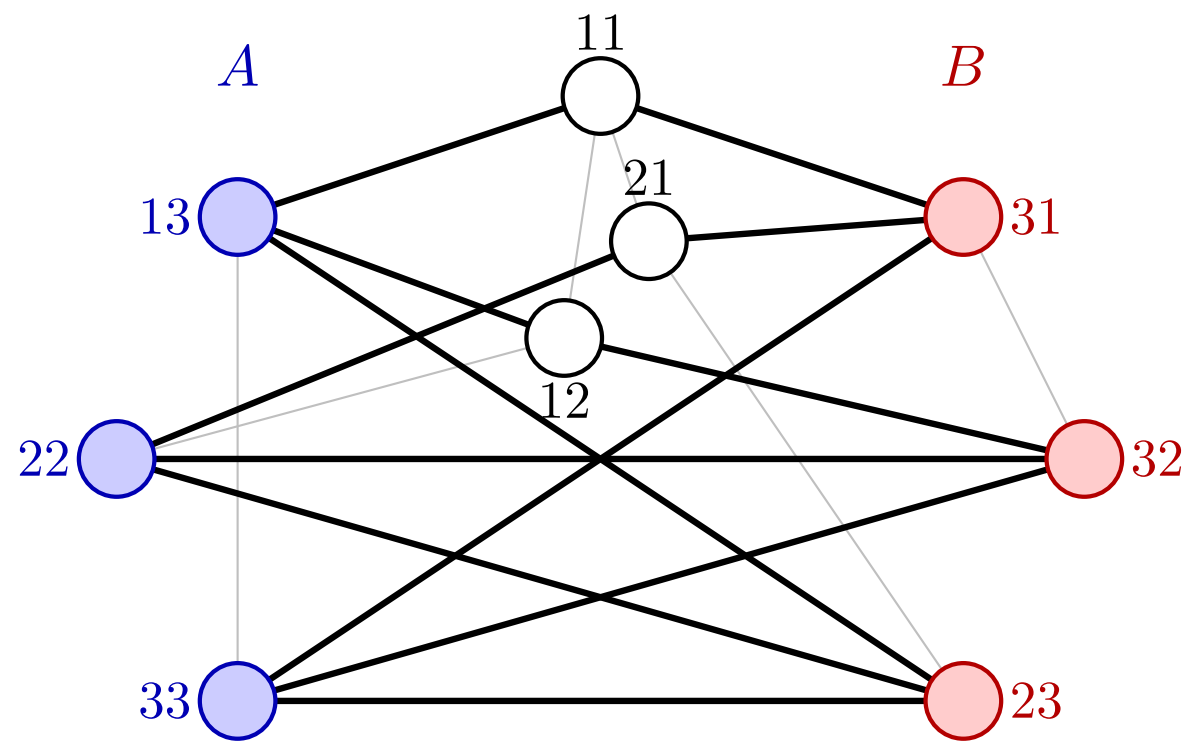
\begin{figure}[ht]
		\centering
		\begin{tikzpicture}[
			scale=0.8,
			branchA/.style={
				circle,
				draw=blue!70!black,
				fill=blue!20,
				thick,
				minimum size=5mm,
				inner sep=0pt
			},
			branchB/.style={
				circle,
				draw=red!70!black,
				fill=red!20,
				thick,
				minimum size=5mm,
				inner sep=0pt
			},
			internal/.style={
				circle,
				draw=black,
				fill=white,
				thick,
				minimum size=5mm,
				inner sep=0pt
			},
			vertexLabel/.style={
				font=\small,
				inner sep=1pt
			},
			graphEdge/.style={
				gray!50,
				thin
			},
			subdivisionEdge/.style={
				black,
				very thick
			}
			]
			
			\coordinate (p13) at (-3, 2);
			\coordinate (p22) at (-4, 0);
			\coordinate (p33) at (-3,-2);
			
			\coordinate (p31) at (3, 2);
			\coordinate (p32) at (4, 0);
			\coordinate (p23) at (3,-2);
			
			\coordinate (p11) at (0,3);
			\coordinate (p12) at (-0.3,1);
			\coordinate (p21) at (0.4,1.8);
			
			\draw[graphEdge] (p11) -- (p12);
			\draw[graphEdge] (p21) -- (p23);
			\draw[graphEdge] (p11) -- (p21);
			\draw[graphEdge] (p12) -- (p22);
			\draw[graphEdge] (p13) -- (p33);
			\draw[graphEdge] (p12) -- (p13);
			\draw[graphEdge] (p31) -- (p32);
			
			\draw[subdivisionEdge] (p13) -- (p11);
			\draw[subdivisionEdge] (p11) -- (p31);
			
			\draw[subdivisionEdge] (p13) -- (p12);
			\draw[subdivisionEdge] (p12) -- (p32);
			
			\draw[subdivisionEdge] (p22) -- (p21);
			\draw[subdivisionEdge] (p21) -- (p31);
			
			\draw[subdivisionEdge] (p13) -- (p23);
			
			\draw[subdivisionEdge] (p22) -- (p32);
			\draw[subdivisionEdge] (p22) -- (p23);
			
			\draw[subdivisionEdge] (p33) -- (p31);
			\draw[subdivisionEdge] (p33) -- (p32);
			\draw[subdivisionEdge] (p33) -- (p23);
			
			\node[
			branchA,
			label={[vertexLabel,blue!70!black]left:$13$}
			] at (p13) {};
			
			\node[
			branchA,
			label={[vertexLabel,blue!70!black]left:$22$}
			] at (p22) {};
			
			\node[
			branchA,
			label={[vertexLabel,blue!70!black]left:$33$}
			] at (p33) {};
			
			\node[
			branchB,
			label={[vertexLabel,red!70!black]right:$31$}
			] at (p31) {};
			
			\node[
			branchB,
			label={[vertexLabel,red!70!black]right:$32$}
			] at (p32) {};
			
			\node[
			branchB,
			label={[vertexLabel,red!70!black]right:$23$}
			] at (p23) {};
			
			\node[
			internal,
			label={[vertexLabel]above:$11$}
			] at (p11) {};
			
			\node[
			internal,
			label={[vertexLabel]below:$12$}
			] at (p12) {};
			
			\node[
			internal,
			label={[vertexLabel]above:$21$}
			] at (p21) {};
			
			\node[
			blue!70!black,
			font=\bfseries
			] at (-3,3.25) {$A$};
			
			\node[
			red!70!black,
			font=\bfseries
			] at (3,3.25) {$B$};
			
		\end{tikzpicture}
		
		\caption{The graph $F^2_{2\times1}(K_3)$. The thick black
			edges form a subdivision of $K_{3,3}$ with branch classes
			$A=\{13,22,33\}$ and $B=\{31,32,23\}$. The white vertices
			$11,12,21$ are the internal vertices of the subdivided paths.}
		\label{fig:F2-2-K3-K33-subdivision}
	\end{figure}



\begin{theorem}
	\label{thm:planarity-complete-supertoken}
	Let $n\geq3$, let $s\geq2$, and let
	$
	2\leq k\leq ns.
	$
	Then $	F^s_{k\times1}(K_n)$ is planar
	if and only if
	$k=ns$.
\end{theorem}

\begin{proof}
	Suppose first that $k=ns$. Every vertex of $K_n$ then contains
	exactly $s$ tokens in every configuration. Thus no token can be
	moved, and $F^s_{ns\times1}(K_n)$ is edgeless. In particular, it is
	planar. Conversely, assume that $2\leq k<ns$. We show that $F^s_{k\times1}(K_n)$ contains one of the nonplanar
	graphs from Lemma~\ref{lem:small-nonplanar-supertoken} as a subgraph.
	
	Choose three distinct vertices $u,v,w\in V(K_n)$ and distinguish a suitable collection of tokens, which we call the
	\emph{active tokens}. All other tokens will be kept fixed. We
	distinguish two cases.

    \medskip
\noindent
\textbf{Case 1: $k\leq ns-4$.} 	Choose two active tokens. We claim that the remaining $k-2$ tokens
	can be placed so that each of $u,v,w$ has at least two vacant places. To obtain such a placement, allow at most $s-2$ fixed tokens at each
	of $u,v,w$ and at most $s$ fixed tokens at each of the other $n-3$ vertices. The number of places available for the fixed tokens
	is 	$3(s-2)+(n-3)s=ns-6$. Since $k\leq ns-4$, we have
	$
	k-2\leq ns-6.
	$
	Hence the number of fixed tokens does not exceed the number of
	places available for them. Consequently, the $k-2$ fixed tokens can be distributed among the
	vertices so that each of $u,v,w$ contains at most $s-2$ fixed tokens,
	while every other vertex contains at most $s$ fixed tokens. 
	Now consider only those configurations in which the fixed tokens
	remain in their chosen positions and each of the two active tokens
	occupies one of $u,v,w$. Since at least two places are available at
	each of these three vertices, the two active tokens may be placed
	arbitrarily on $u,v,w$, including at the same vertex. As $K_n$ is complete, moving either active token from one of
	$u,v,w$ to another is always legal. Hence, the subgraph induced by
	these configurations is isomorphic to $F^2_{2\times1}(K_3)$. By Lemma~\ref{lem:small-nonplanar-supertoken}, this subgraph is
	nonplanar. Consequently, $F^s_{k\times1}(K_n)$ is nonplanar.

    \medskip
\noindent
\textbf{Case 2: $ns-3\leq k<ns$.} Put $c:=ns-k$. 	The integer $c$ is the number of unoccupied places in any
	configuration, because the graph has total capacity $ns$ and
	contains $k$ tokens. From $ns-3\leq k<ns$, we obtain
	$
	1\leq ns-k\leq3,
	$
	and hence
	$
	c\in\{1,2,3\}.
	$ At each of $u,v,w$, place $s-2$ fixed tokens. At each of the
	remaining $n-3$ vertices, place $s$ fixed tokens. The number of fixed
	tokens is $3(s-2)+(n-3)s=ns-6$. Therefore, the number of tokens not yet placed is $
	q
	=
	k-(ns-6)
	=
	6-c.
	$
	Since $c\in\{1,2,3\}$, it follows that
	$
	q\in\{3,4,5\}.
	$
	Declare the remaining $q$ tokens to be the \emph{active tokens}, and
	consider only those configurations in which all fixed tokens remain
	in their prescribed positions, and every active token occupies one of
	the vertices $u,v,w$. Each vertex outside $\{u,v,w\}$ contains exactly $s$ fixed tokens
	and is therefore full. Hence no active token can be moved to any
	vertex outside $\{u,v,w\}$. At each of $u,v,w$, however, we placed
	exactly $s-2$ fixed tokens. Consequently, each of these three
	vertices has precisely two places available for the active tokens.
	Thus, an active configuration is admissible if and only if no more
	than two active tokens occupy any one of $u,v,w$. Furthermore, since $K_n$ is complete, the subgraph induced by
	$\{u,v,w\}$ is isomorphic to $K_3$. Therefore, an active token can move
	from any one of $u,v,w$ to either of the other two, provided that the
	destination currently contains fewer than two active tokens. 
    These are precisely the legal moves for $q$ distinguishable tokens on $K_3$, with a capacity of $2$ tokens per vertex. It follows that, after omitting the fixed tokens from the notation,
	the selected configurations and their adjacencies form a graph
	isomorphic to $F^2_{q\times1}(K_3)$ for $q\in\{3,4,5\}$.
	By Lemma~\ref{lem:small-nonplanar-supertoken}, each of these graphs is nonplanar. Hence $F^s_{k\times1}(K_n)$ contains a nonplanar
	subgraph and is itself nonplanar.
	\end{proof}

\subsection{Classification for complete graphs}

    We conclude this section by giving a complete classification for all
capacities. 
The case $s=1$ was settled by Dalf\'o, Fern\'andez, Fiol, and Reyes~\cite{dffr26}. 


\begin{theorem}[{\cite[Theorem~10.5]{dffr26}}]
	\label{thm:complete-capacity-one-planarity}
	Let \(n\geq2\) and \(1\leq k\leq n-1\). Then
	\(F^{1}_{k\times1}(K_n)\) is planar if and only if
	$
	(n,k)\in
	\{(2,1),(3,1),(3,2),(4,1),(4,2)\}.
	$
\end{theorem}

Combining the preceding theorem with our results for \(s\geq2\) completes the planarity classification of supertoken graphs over complete graphs. 
The remaining planar cases are those listed below, where
planarity occurs only at full capacity, when every vertex is saturated, and the supertoken graph is edgeless.

\begin{theorem}
	\label{thm:complete-planarity-classification}
	Let \(n,s\geq1\) and \(1\leq k\leq ns\). Then
	\(F^s_{k\times1}(K_n)\) is planar if and only if one of the following
	conditions holds:
	\begin{enumerate}[label=\textnormal{(\roman*)}]
		\item \(n=1\);
		
		\item \(n\geq2\), \(s=1\), and either \(k=n\) or
		\((n,k)\in\{(2,1),(3,1),(3,2),(4,1),(4,2)\}\);
		
		\item \(n=2\), \(s\geq2\), and either
		\(k\in\{1,2,3,2s\}\) or \((s,k)=(3,4)\);
		
		\item \(n\geq3\), \(s\geq2\), and either \(k=ns\), or
		\(k=1\) and \(n\leq4\).
	\end{enumerate}
\end{theorem}

\begin{proof}
	If \(n=1\), then \(F^s_{k\times1}(K_1)\) consists of a single
	vertex, and hence is planar. Assume \(n\geq2\). If \(s=1\), then for \(k<n\) the result follows
	directly from Theorem~\ref{thm:complete-capacity-one-planarity},
	whereas for \(k=n\) the graph is edgeless. This gives
	condition~\textnormal{(ii)}. Now let \(s\geq2\). If \(k=1\), then
	\(F^s_{1\times1}(K_n)\cong K_n\), which is planar if and only if
	\(n\leq4\). For \(k\geq2\), if \(n=2\), condition~\textnormal{(iii)} follows
	from Theorem~\ref{thm:planarity-supertoken-K2}, together with the
	case \(k=1\) above. If \(n\geq3\), condition~\textnormal{(iv)}
	follows from Theorem~\ref{thm:planarity-complete-supertoken},
	again together with the case \(k=1\). Thus the four conditions above exhaust all possibilities.
\end{proof}

\subsection*{Acknowledgments} 
The first author is supported by the China Scholarship Council program (Project ID: 202606750007).  
The second and fourth authors' research was supported by project PID2025-171721NB-I00 funded by MICIU/AEI/10.13039/501100011033 and by ERDF, EU; and by project
2021SGR00434 funded by AGAUR. 
In addition, the third author's research was partially supported by the Slovenian Research and Innovation Agency (ARIS) under research program P1-0285 and research projects J1-50000 and Z1-70008.

\subsection*{Conflict of interest}
The authors declare no conflict of interest. 

\subsection*{Data availability statement} 
All data generated or analyzed during this study are included in this article.



\section*{Appendix}
\label{app:pseudocode}

\begin{samepage}
We give below the pseudo-code used to compute exactly the diameter of
$F^s_{k\times1}(C_n)$ for odd $n$. The algorithm uses the symmetries
of $C_n$ and breadth-first search to reduce the computation.

\begin{algorithm}[H]
\caption{Diameter of $F^s_{k\times 1}(G)$}
\label{alg:diameter}
\begin{algorithmic}[1]
\Require A graph $G=(V,E)$ and integers $s,k$ with $1\leq k\leq s|V|$
\Ensure $\operatorname{diam}(F^s_{k\times1}(G))$

\State $\mathcal{X}\gets
\{(x_1,\ldots,x_k)\in V^k:
|\{i:x_i=v\}|\leq s\text{ for every }v\in V\}$

\State Choose a set $\mathcal{R}\subseteq\mathcal{X}$ containing one
representative of each orbit under token permutations and automorphisms of $G$
\State $D\gets 0$

\For{$X\in\mathcal{R}$}
    \State $d(X)\gets 0$ and $d(Y)\gets\infty$ for every $Y\in\mathcal{X}\setminus\{X\}$
    \State initialize a queue $Q$ with $X$
    \While{$Q\neq\emptyset$}
        \State remove the first configuration $Y=(y_1,\ldots,y_k)$ from $Q$
        \For{$i=1,\ldots,k$}
            \For{$v\in N_G(y_i)$ such that $|\{j:y_j=v\}|<s$}
                \State $Z\gets(y_1,\ldots,y_{i-1},v,y_{i+1},\ldots,y_k)$
                \If{$d(Z)=\infty$}
                    \State $d(Z)\gets d(Y)+1$; append ${\boldsymbol{z}}$ to $Q$
                \EndIf
            \EndFor
        \EndFor
    \EndWhile
    \If{$d(Y)=\infty$ for some $Y\in\mathcal{X}$}
        \State \Return $+\infty$
    \EndIf
    \State $D\gets\max\{D,\max_{Y\in\mathcal{X}}d(Y)\}$
\EndFor

\State \Return $D$
\end{algorithmic}
\end{algorithm}
\end{samepage}

\subsection*{Subdivisions of $K_{3,3}$ in $F^2_{q\times1}(K_3)$}

To support the proof of Lemma~\ref{lem:small-nonplanar-supertoken}, we
explicitly exhibit a subdivision of $K_{3,3}$ in
$F^2_{q\times1}(K_3)$ for each $q\in\{3,4,5\}$.

In each of the figures below, the blue and red vertices form the two
branch sets $A$ and $B$, respectively, each of cardinality three. The
thick black paths join every vertex of $A$ to every vertex of $B$, and
their internal vertices, shown in white, are pairwise distinct and
disjoint from $A\cup B$. Hence, the union of these nine paths forms a
subdivision of $K_{3,3}$.

Moreover, every configuration appearing in the figures is admissible,
and consecutive configurations along each path differ in exactly one
coordinate. Since the base graph is $K_3$, each such change corresponds
to a permissible token move. Thus all the displayed paths are contained
in the corresponding supertoken graph.

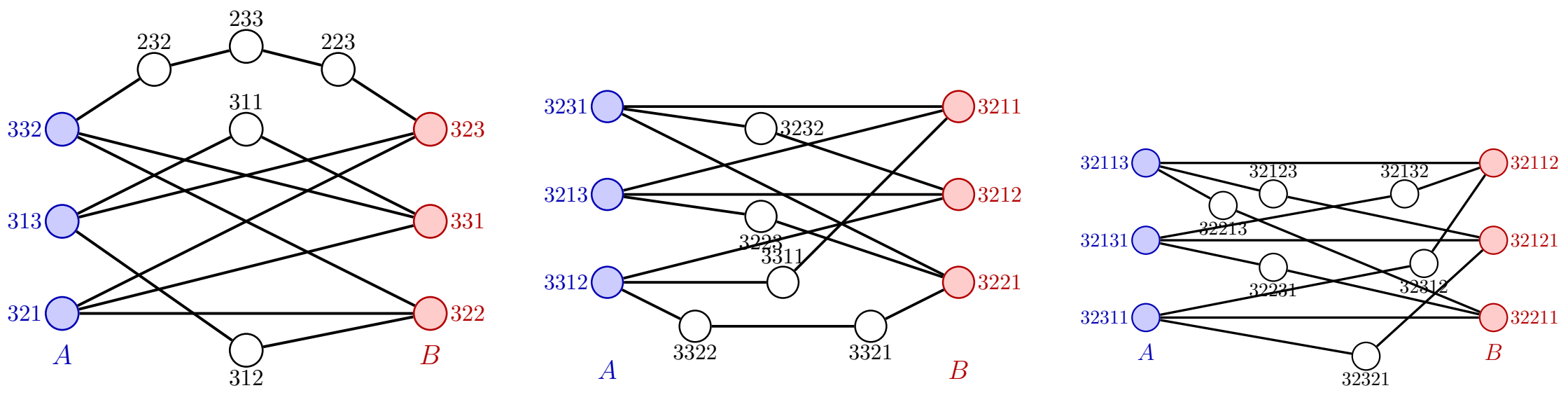
\begin{figure}[ht]
\centering

\begin{subfigure}[t]{0.31\textwidth}
\centering
\resizebox{\linewidth}{!}{%
\begin{tikzpicture}[
    scale=0.7,
    branchA/.style={
        circle,
        draw=blue!70!black,
        fill=blue!20,
        thick,
        minimum size=5mm,
        inner sep=0pt
    },
    branchB/.style={
        circle,
        draw=red!70!black,
        fill=red!20,
        thick,
        minimum size=5mm,
        inner sep=0pt
    },
    internal/.style={
        circle,
        draw=black,
        fill=white,
        thick,
        minimum size=5mm,
        inner sep=0pt
    },
    vertexLabel/.style={
        font=\small,
        inner sep=1pt
    },
    subdivisionEdge/.style={
        black,
        very thick
    }
]

\node[
branchA,
label={[vertexLabel,blue!70!black]left:$332$}
] (a1) at (0,4) {};

\node[
branchA,
label={[vertexLabel,blue!70!black]left:$313$}
] (a2) at (0,2) {};

\node[
branchA,
label={[vertexLabel,blue!70!black]left:$321$}
] (a3) at (0,0) {};

\node[
branchB,
label={[vertexLabel,red!70!black]right:$323$}
] (b1) at (8,4) {};

\node[
branchB,
label={[vertexLabel,red!70!black]right:$331$}
] (b2) at (8,2) {};

\node[
branchB,
label={[vertexLabel,red!70!black]right:$322$}
] (b3) at (8,0) {};

\node[
internal,
label={[vertexLabel]above:$232$}
] (p1) at (2,5.3) {};

\node[
internal,
label={[vertexLabel]above:$233$}
] (p2) at (4,5.8) {};

\node[
internal,
label={[vertexLabel]above:$223$}
] (p3) at (6,5.3) {};

\node[
internal,
label={[vertexLabel]above:$311$}
] (p4) at (4,4) {};

\node[
internal,
label={[vertexLabel]below:$312$}
] (p5) at (4,-0.8) {};

\draw[subdivisionEdge]
(a1) -- (p1) -- (p2) -- (p3) -- (b1);

\draw[subdivisionEdge]
(a2) -- (p4) -- (b2);

\draw[subdivisionEdge]
(a2) -- (p5) -- (b3);

\draw[subdivisionEdge] (a1) -- (b2);
\draw[subdivisionEdge] (a1) -- (b3);
\draw[subdivisionEdge] (a2) -- (b1);

\draw[subdivisionEdge] (a3) -- (b1);
\draw[subdivisionEdge] (a3) -- (b2);
\draw[subdivisionEdge] (a3) -- (b3);

\node[
blue!70!black,
font=\large\bfseries
] at (0,-0.9) {$A$};

\node[
red!70!black,
font=\large\bfseries
] at (8,-0.9) {$B$};

\end{tikzpicture}%
}
\caption{$q=3$.}
\label{fig:subdivision-K33-q3}
\end{subfigure}
\hfill
\begin{subfigure}[t]{0.31\textwidth}
\centering
\resizebox{\linewidth}{!}{%
\begin{tikzpicture}[
    scale=0.7,
    branchA/.style={
        circle,
        draw=blue!70!black,
        fill=blue!20,
        thick,
        minimum size=5mm,
        inner sep=0pt
    },
    branchB/.style={
        circle,
        draw=red!70!black,
        fill=red!20,
        thick,
        minimum size=5mm,
        inner sep=0pt
    },
    internal/.style={
        circle,
        draw=black,
        fill=white,
        thick,
        minimum size=5mm,
        inner sep=0pt
    },
    vertexLabel/.style={
        font=\small,
        inner sep=1pt
    },
    subdivisionEdge/.style={
        black,
        very thick
    }
]

\node[
branchA,
label={[vertexLabel,blue!70!black]left:$3231$}
] (a1) at (0,4) {};

\node[
branchA,
label={[vertexLabel,blue!70!black]left:$3213$}
] (a2) at (0,2) {};

\node[
branchA,
label={[vertexLabel,blue!70!black]left:$3312$}
] (a3) at (0,0) {};

\node[
branchB,
label={[vertexLabel,red!70!black]right:$3211$}
] (b1) at (8,4) {};

\node[
branchB,
label={[vertexLabel,red!70!black]right:$3212$}
] (b2) at (8,2) {};

\node[
branchB,
label={[vertexLabel,red!70!black]right:$3221$}
] (b3) at (8,0) {};

\node[
internal,
label={[vertexLabel]right:$3232$}
] (p1) at (3.5,3.5) {};

\node[
internal,
label={[vertexLabel]below:$3223$}
] (p2) at (3.5,1.5) {};

\node[
internal,
label={[vertexLabel]above:$3311$}
] (p3) at (4,0) {};

\node[
internal,
label={[vertexLabel]below:$3322$}
] (p4) at (2,-1) {};

\node[
internal,
label={[vertexLabel]below:$3321$}
] (p5) at (6,-1) {};

\draw[subdivisionEdge]
(a1) -- (p1) -- (b2);

\draw[subdivisionEdge]
(a2) -- (p2) -- (b3);

\draw[subdivisionEdge]
(a3) -- (p3) -- (b1);

\draw[subdivisionEdge]
(a3) -- (p4) -- (p5) -- (b3);

\draw[subdivisionEdge] (a1) -- (b1);
\draw[subdivisionEdge] (a1) -- (b3);

\draw[subdivisionEdge] (a2) -- (b1);
\draw[subdivisionEdge] (a2) -- (b2);

\draw[subdivisionEdge] (a3) -- (b2);

\node[
blue!70!black,
font=\large\bfseries
] at (0,-2) {$A$};

\node[
red!70!black,
font=\large\bfseries
] at (8,-2) {$B$};

\end{tikzpicture}%
}
\caption{$q=4$.}
\label{fig:subdivision-K33-q4}
\end{subfigure}
\hfill
\begin{subfigure}[t]{0.31\textwidth}
\centering
\resizebox{\linewidth}{!}{%
\begin{tikzpicture}[
    scale=0.7,
    branchA/.style={
        circle,
        draw=blue!70!black,
        fill=blue!20,
        thick,
        minimum size=5mm,
        inner sep=0pt
    },
    branchB/.style={
        circle,
        draw=red!70!black,
        fill=red!20,
        thick,
        minimum size=5mm,
        inner sep=0pt
    },
    internal/.style={
        circle,
        draw=black,
        fill=white,
        thick,
        minimum size=5mm,
        inner sep=0pt
    },
    vertexLabel/.style={
        font=\small,
        inner sep=1pt
    },
    subdivisionEdge/.style={
        black,
        very thick
    }
]

\node[
branchA,
label={[vertexLabel,blue!70!black]left:$32113$}
] (a1) at (0,4) {};

\node[
branchA,
label={[vertexLabel,blue!70!black]left:$32131$}
] (a2) at (0,2) {};

\node[
branchA,
label={[vertexLabel,blue!70!black]left:$32311$}
] (a3) at (0,0) {};

\node[
branchB,
label={[vertexLabel,red!70!black]right:$32112$}
] (b1) at (9,4) {};

\node[
branchB,
label={[vertexLabel,red!70!black]right:$32121$}
] (b2) at (9,2) {};

\node[
branchB,
label={[vertexLabel,red!70!black]right:$32211$}
] (b3) at (9,0) {};

\node[
internal,
label={[vertexLabel]above:$32123$}
] (p12) at (3.3,3.2) {};

\node[
internal,
label={[vertexLabel]below:$32213$}
] (p13) at (2,2.9) {};

\node[
internal,
label={[vertexLabel]above:$32132$}
] (p21) at (6.7,3.2) {};

\node[
internal,
label={[vertexLabel]below:$32231$}
] (p23) at (3.3,1.3) {};

\node[
internal,
label={[vertexLabel]below:$32312$}
] (p31) at (7.2,1.4) {};

\node[
internal,
label={[vertexLabel]below:$32321$}
] (p32) at (5.7,-1) {};

\draw[subdivisionEdge]
(a1) -- (p12) -- (b2);

\draw[subdivisionEdge]
(a1) -- (p13) -- (b3);

\draw[subdivisionEdge]
(a2) -- (p21) -- (b1);

\draw[subdivisionEdge]
(a2) -- (p23) -- (b3);

\draw[subdivisionEdge]
(a3) -- (p31) -- (b1);

\draw[subdivisionEdge]
(a3) -- (p32) -- (b2);

\draw[subdivisionEdge] (a1) -- (b1);
\draw[subdivisionEdge] (a2) -- (b2);
\draw[subdivisionEdge] (a3) -- (b3);

\node[
blue!70!black,
font=\large\bfseries
] at (0,-0.9) {$A$};

\node[
red!70!black,
font=\large\bfseries
] at (9,-0.9) {$B$};

\end{tikzpicture}%
}
\caption{$q=5$.}
\label{fig:subdivision-K33-q5}
\end{subfigure}

\caption{Subdivisions of $K_{3,3}$ contained in
$F^2_{q\times1}(K_3)$ for $q=3,4,5$. In each panel, the blue and red
vertices form the branch sets $A$ and $B$, respectively, and the white
vertices are the internal vertices of the subdivided edges.}
\label{fig:subdivision-K33}
\end{figure}

The subdivisions, for $q=3,4,5$, are shown in
Figures~\ref{fig:subdivision-K33}.

\end{document}